\documentclass[12pt,a4paper,reqno]{amsart}

\usepackage[english]{babel}
\usepackage{amssymb,latexsym,amsfonts,amsthm,upref,amsmath}
\usepackage[foot]{amsaddr}
\usepackage[margin=1in]{geometry} 
\usepackage[dvipsnames,x11names]{xcolor}
 \definecolor{myblue}{HTML}{003399}
\usepackage{hyperref}
\hypersetup{colorlinks,citecolor=myblue,filecolor=black,linkcolor=myblue,urlcolor=myblue}
\usepackage{enumerate}  
\usepackage{tikz}
\usepackage{float}
\usepackage{cleveref}
\usepackage{mathtools}

\DeclarePairedDelimiter\floor{\lfloor}{\rfloor}
\usepackage{cite}
\usepackage{filecontents}
\makeatletter
\newcommand{\leqnomode}{\tagsleft@true}
\newcommand{\reqnomode}{\tagsleft@false}

\makeatother
\newtheorem*{thm*}{Theorem}
\newtheorem*{lem*}{Lemma}
\newtheoremstyle{prim}{}{}{\normalfont}{}{\bfseries}{.}{ }{}
\newtheoremstyle{stil}{}{}{\slshape}{}{\bfseries}{.}{ }{}
\theoremstyle{stil}
\newtheorem{thm}{Theorem}[section]
\newtheoremstyle{defi}{}{}{}{}{\bfseries}{.}{ }{}
\theoremstyle{defi}
\newtheorem{defn}[thm]{Definition}
\theoremstyle{defi}
\newtheorem{rem}[thm]{Remark}
\theoremstyle{stil}
\newtheorem*{mthm*}{Main Theorem}
\newtheorem*{kor*}{Corollary}
\newtheorem{pro}[thm]{Proposition}
\theoremstyle{stil}
\newtheorem{lem}[thm]{Lemma}
\theoremstyle{stil}
\newtheorem{kor}[thm]{Corollary}
\theoremstyle{prim}
\newtheorem{ex}[thm]{Example}
\newenvironment{prf}{\noindent \textit{Proof.}}{\null\hfill$\qed$\hskip
2mm\vskip 2mm}

\newcommand{\wa}{ \hat{a}}
\newcommand{\wy}{ \hat{y}}
\newcommand{\wx}{ \hat{x}}
\newcommand{\hp}{ {\rm H}^+_N}

\newcommand{\hh}{ {\rm H}}

\newcommand{\vac}{ \mathrm{\boldsymbol{1}}}

\newcommand{\h}{\mathfrak h}
\newcommand{\hhat}{\hat{\mathfrak{h}}}

\newcommand{\CC}{\mathbb{C}}

\newcommand{\ZZ}{\mathbb{Z}}

\newcommand{\Pc}{\mathcal{P}}

\newcommand{\Sc}{\mathcal{S}}
\newcommand{\Sd}{\wvr{S}}

\newcommand{\Vc}{\mathcal{V}}

\newcommand{\wvr}{\overline}

\newcommand{\ot}{\otimes}
\newcommand{\ts}{\hspace{1pt}}

\newcommand{\tr}{ {\rm tr}}

\newcommand{\xpan}{\mathop{\mathrm{span}}}
\newcommand{\ndo}{\mathop{\mathrm{End}}}
\newcommand{\om}{\mathop{\mathrm{Hom}}}

\newcommand{\rez}{\mathop{\mathrm{Res}}}

\newcommand{\diag}{\mathop{\mathrm{diag}}}

\newcommand{\fand}{\quad\text{and}\quad}
\newcommand{\Fand}{\qquad\text{and}\qquad}

\newcommand{\non}{\nonumber}
\newcommand{\beq}{\begin{equation}}
\newcommand{\eeq}{\end{equation}}
\newcommand{\ben}{\begin{equation*}}
\newcommand{\een}{\end{equation*}}

\makeatletter
\def\smalloverbrace#1{\mathop{\vbox{\m@th\ialign{##\crcr\noalign{\kern3\p@}%
  \tiny\downbracefill\crcr\noalign{\kern3\p@\nointerlineskip}%
  $\hfil\displaystyle{#1}\hfil$\crcr}}}\limits}
\makeatother

\makeatletter
\def\smallunderbrace#1{\mathop{\vtop{\m@th\ialign{##\crcr
   $\hfil\displaystyle{#1}\hfil$\crcr
   \noalign{\kern3\p@\nointerlineskip}%
   \tiny\upbracefill\crcr\noalign{\kern3\p@}}}}\limits}
\makeatother

\makeatletter
\def\author@andify{%
  \nxandlist {\unskip ,\penalty-1 \space\ignorespaces}%
    {\unskip {} \@@and~}%
    {\unskip \penalty-2 \space \@@and~}%
}
\makeatother

\begin{document}

\title{On the Heisenberg algebra associated with the rational \texorpdfstring{$R$}{R}-matrix}

\author{Marijana Butorac}
\address[M. Butorac]{Department of Mathematics, University of Rijeka, Radmile Matej\v{c}i\'{c} 2, 51\,000 Rijeka, Croatia}
\email{mbutorac@math.uniri.hr}

\author{Slaven Ko\v{z}i\'{c}}
\address[S. Ko\v{z}i\'{c}]{Department of Mathematics, Faculty of Science, University of Zagreb,  Bijeni\v{c}ka cesta 30, 10\,000 Zagreb, Croatia}
\email{kslaven@math.hr}

\begin{abstract}
We associate a deformation of Heisenberg algebra   to the suitably normalized Yang $R$-matrix  and we investigate its properties. Moreover, we construct   new examples   of quantum vertex algebras which possess  the same representation theory
as the aforementioned deformed Heisenberg algebra.
\end{abstract}

\maketitle

\allowdisplaybreaks


\section{Introduction}\label{intro}
\numberwithin{equation}{section}

The notion of vertex algebra, which was  introduced by Borcherds \cite{B},  presents a remarkable connection between  mathematics and theoretical physics. Starting with   Belavin, Polyakov and Zamolodchikov \cite{BPZ}, such objects were extensively studied by physicists in connection with  conformal symmetries of two-dimensional quantum field theory. On the other hand, the theory of vertex algebras led to important new  methods, techniques and results in multiple areas of mathematics such as affine Kac--Moody Lie algebras, automorphic forms, finite simple groups and $\mathcal{W}$-algebras; see, e.g., the books by E.\,Frenkel and Ben-Zvi \cite{FBZ}, I. Frenkel, Lepowsky and Meurman \cite{FLM} and Kac \cite{Kac2}.

Let $\h$ be an  abelian Lie algebra over $\CC$ equipped with the  nondegenerate invariant symmetric
bilinear form $\left<\cdot,\cdot\right>$.
The   affine Lie algebra $\hhat$ is defined on the complex   space
\beq\label{introabafa}
\hhat=\h\ot\CC[t,t^{-1}] \oplus\CC C 
\eeq
with the Lie brackets   given by
\beq\label{introbrckts}
\left[a(r), b(s)\right] =\left<a,b\right> r\ts\delta_{r+s\ts 0}\ts C
\fand
\left[C,x\right]=0 
\eeq
for all $r,s\in\ZZ$, $a,b\in\h$ and $x\in\hhat$, where $a(r)$  denotes the element $a\ot t^r$. The corresponding Heisenberg Lie algebra $\hhat_*$ is defined as a  subalgebra
$$
\hhat_* =\coprod_{n\in\ZZ\setminus\left\{0\right\}} \left(\h\ot t^n\right)\oplus\CC C \subset \hhat.
$$
As with many other infinite-dimensional Lie algebras,   Heisenberg Lie algebras  play   an important role in the theory of vertex algebras and their representations; see, e.g.,  \cite{FK,FZ,LW,S}.
In this paper, we study   certain deformation  of the universal enveloping algebra of the Heisenberg Lie algebra, along with the underlying  quantum  vertex algebra theory. 

The notion of quantum vertex algebra was introduced by Etingof and Kazhdan
in \cite{EK},  where they also  associated the  quantum affine vertex algebras   to the rational, trigonometric and elliptic $R$-matrix of type $A$. The  $\Sc$-locality of these quantum   vertex algebras, i.e. the quantum version of the locality property for vertex algebras,   possesses the form of the so-called quantum current commutation relation, which   goes back to Reshetikhin and  Semenov-Tian-Shansky \cite{RS}. 
In this paper, we continue the study \cite{c11} of the interplay between the quantum current commutation relation  associated with the rational $R$-matrix of type $A$  and quantum vertex algebra theory.
Motivated by the form of the Heisenberg Lie algebra defining relation \eqref{introbrckts}, we investigate quantum (vertex) algebras   defined by the   relations which come from the aforementioned quantum current commutation  relation  by extracting   the quadratic terms and   terms containing only the central element $C$ or the unit $1$. 
Such relations can be expressed in the form similar to  \eqref{introbrckts} as
\beq\label{introrel}
y_1(u) \ts y_2 (v) +S_{12} (u-v,C)=y_2 (v)\ts y_1(u) +S_{21} (v-u,C),
\eeq
where $y(z)$ denotes the $N\times N$ matrix of formal power series of the algebra generators and $S(z,C)$ is a certain product of two copies of the suitably normalized Yang $R$-matrices. The precise meaning of \eqref{introrel} is explained in Subsection \ref{subsec20}.

In Section \ref{sec2},  we construct    quantum vertex algebras $\hp$ whose braiding map $\Sc=\Sc(z)$ is governed by $S(z,c)$, where $c\in\CC$, so that the form of their $\Sc$-locality property resembles \eqref{introrel}.
Moreover, we show that $\hp$ contains a quantum vertex subalgebra $\Vc_{\hh}(c)$   whose classical limit  $h\to 0$ coincides with the level $c$ Heisenberg vertex algebra associated to $\hhat_*$, where $\h$ is a Cartan subalgebra of $\mathfrak{sl}_N$.

 In Section \ref{sec3}, we study certain associative algebras  $\hh(C)$ and $\hh(C)_* $, defined over $\CC[[h]]$, whose defining relations are found by taking the diagonal entries of \eqref{introrel}. In particular, we establish the Poincar\'{e}--Birkhoff--Witt theorem for these algebras. Roughly speaking, $\hh(C)$ and $\hh(C)_*$ can be regarded as  deformations of universal enveloping algebras $U(\hhat)$ and $U(\hhat_*)$, respectively. Next, we turn to  their representation theory and, following the classical theory, we introduce the notion of restricted module. Furthermore, we construct examples of such modules by generalizing the well-known canonical realization of the affine Lie algebra $\hhat$. Finally, we show that the (irreducible) modules for the quantum vertex algebra $\Vc_{\hh}(c)$ coincide with the level $c$ (irreducible) restricted $\hh(C)$-modules.

In the end, we should mention that the problem of associating quantum vertex algebras to certain deformed Heisenberg Lie algebras, which differ from those considered in this paper, was   studied by Li \cite{Li-h}.

\section{Quantum vertex algebras}\label{sec2}
 
In Subsections \ref{subsec20} and \ref{subsec22}, we introduce the  data which is required to define a structure of quantum vertex algebra  over a certain quotient $\hp$ of the $h$-adically completed algebra of polynomials in infinitely many variables. In Subsection \ref{subsec23}, we present the main result of this section, i.e. the aforementioned construction of quantum vertex algebra $\hp$, and in Subsection \ref{subsec24} we give its proof.
Finally, in Subsection \ref{subsec21}, we discuss certain quantum vertex subalgebra $\Vc_{\hh}(c)\subset\hp$ which is a deformation of the Heisenberg vertex algebra.

\subsection{Creation and annihilation operators}\label{subsec20}

Let $N\geqslant 2$ be an integer and
 $h$   a formal parameter. The {\em Yang $R$-matrix} $R(u)=R_{12}(u) \in\ndo\CC^N\ot\ndo\CC^N [h/u]$ is defined by
$$
R(u)=I-\textstyle \frac{h}{u}P,
$$
where $I$ is the identity and $P$ the permutation operator,
\beq\label{perm8}
I=\sum_{i,j=1}^N e_{ii}\ot e_{jj}\fand
P=\sum_{i,j=1}^N e_{ij}\ot e_{ji},
\eeq
and $e_{ij}$ are the matrix units.
Let $C$ be another formal parameter. There exists a unique formal power series 
$$
G(u,C)=1+\frac{C+N}{N}\frac{h^2}{u^2}-\frac{C(C+N)}{N}\frac{h^3}{u^3}+\ldots \,\in\CC[C][[h/u]] 
$$
 such that 
\beq\label{tr34}
\tr_1\ts \left(G(u,C) R(u)R(-u-hC)-I\right)
=
\tr_2 \ts  \left(G(u,C) R(u)R(-u-hC)-I\right)
=0,
\eeq
where $\tr_i$ denotes the trace taken over the $i$-th tensor factor. 
Fix $c\in\CC$ and define a power series $S(u,c)\in u^{-2}\ndo\CC^N \ot\ndo\CC^N [[h/u]]$ by
\beq\label{sofu}
S(u,c)=h^{-2}\left(G(u,c)R(u)R(-u-hc)-I\right).
\eeq
We shall often omit the second argument $c$   and write $S(u)$ instead of $S(u,c)$.
One easily checks that $S(u)$ is well-defined, i.e. that the expression $G(u,c)R(u)R(-u-hc)-I$ possesses a zero of order two at $h=0$, so that \eqref{sofu} does not contain any negative powers of the parameter $h$. Moreover, we have
\beq\label{esform}
S(u)\in\frac{c}{N u^2}\left(I-NP\right) + \frac{h}{u^3} \ndo\CC^N \ot\ndo\CC^N  [[h/u]].
\eeq

Consider  the $h$-adically completed algebra of polynomials in variables
$x_{ij}^{(-r)}$, 
$$
 \Pc = \CC[x_{ij}^{(-r)}\,:\, i,j=1,\ldots ,N,\, r\geqslant 1][[h]].
$$
Let $\mathcal{I}$ be the  $h$-adically complete ideal in $\Pc$ generated by the elements
\beq\label{34els}
x_{11}^{(-r)}+x_{22}^{(-r)}+\ldots + x_{NN}^{(-r)}\quad\text{with}\quad r=1,2,\ldots .
\eeq
We define   the   associative algebra $\hp$  over the ring $\CC[[h]]$  as the quotient
$$
\hp = \Pc / \mathcal{I}.
$$
It is clear that $\hp$ is  topologically free, i.e. a torsion-free, separated and $h$-adically complete $\CC[[h]]$-module; see, e.g.,  \cite[Ch. XVI]{Kas} for more information on $\CC[[h]]$-modules.

It will be convenient to arrange the  elements $x_{ij}^{(-r)} $    into matrices of formal power series,
$$
x^+(u)=\sum_{i,j=1}^N e_{ij} \ot x^+_{ij}(u)\in\ndo\CC^N \ot \hp[[u]],\quad\text{where}\quad
x^+_{ij}(u)=\sum_{r\geqslant 1} x_{ij}^{(-r)} u^{r-1}.
$$
 Also, generalizing the above formula, for any integer $n\geqslant 1$ we write 
\beq\label{xn}
x^+_{[n]}(u)=x^+_{1}(u_1)\ldots x^+_{n}(u_n)\fand
x^+_{[n]}(z+u)=x^+_{1}(z+u_1)\ldots x^+_{n}(z+u_n),
\eeq
 where $u=(u_1,\ldots, u_n)$ is a family of variables,  $z$ a  single variable and
\beq\label{notation788}
x^+_{k}(u)=\sum_{i,j=1}^N 1^{\ot (k-1)}\ot e_{ij} \ot 1^{\ot (n-k)}\ot x^+_{ij}(u)
\quad\text{with}\quad k=1,\ldots ,n.
\eeq
Hence the coefficients of the expressions in \eqref{xn} belong to $(\ndo\CC^N)^{\ot n} \ot \hp$. Throughout the paper we shall often use the notation as in \eqref{notation788}, where the subscripts indicate   factors in the tensor product algebra.

 Let $V$ be a $\CC[[h]]$-module. We denote  by $V((u))_h$ the $\CC[[h]]$-module of all series  
\beq\label{notacija}
a(u)=\sum_{r\in\ZZ } a_r u^{-r-1}\in V[[u^{\pm 1}]]\quad \text{such that}\quad a_r \rightarrow 0 \text{ when } r \rightarrow \infty
\eeq
 with respect to the $h$-adic topology. 
Furthermore, we denote by $V[u^{-1}]_h$ the $\CC[[h]]$-module of all series as in \eqref{notacija} such that, in addition, $a(u)$ belongs to $V[[u^{-1}]]$ .
Such notation naturally extends to the multiple variable case,   so that we  write, e.g., $V((u_1,\ldots ,u_n))_h $.
Observe that if   $V$ is a topologically free $\CC[[h]]$-module, hence isomorphic to $V_0[[h]]$ for some complex   space $V_0$, then $V((u))_h$ is  topologically free as well. Moreover, $V((u))_h$ can  be then identified with $V_0((u))[[h]]$, which is   the $h$-adic completion of $V((u))$.

 The next proposition  can be easily proved by  using    \eqref{tr34}  and the defining relations for the algebra $\hp$, 
\beq\label{tr23}
x_{11}^{(-r)}+x_{22}^{(-r)}+\ldots + x_{NN}^{(-r)}=0\quad\text{for all}\quad r=1,2,\ldots .
\eeq
\begin{lem}\label{lemma21}
For any $c\in\CC$ there exists a unique operator  
$$
x^-(u)=\sum_{i,j=1}^N e_{ij} \ot x^-_{ij}(u), \quad\text{where}\quad
x^-_{ij}(u)=\sum_{r\geqslant 1} x_{ij}^{(r-1)} u^{-r},
$$
which belongs to  $\ndo\CC^N \ot \om (\hp,\hp[u^{-1}]_h)$,
such that $x^-(u)\vac=0$ and for any integer  $n\geqslant 1$ and the  variables $u=(u_1,\ldots, u_n)$ we have
\beq\label{xm}
x^-_0(u_0)\ts x^+_{[n]}(u)
=-\sum_{j=1}^n S_{0j}(u_0-u_j) \ts x^+_{1}(u_1)\ldots x^+_{j-1}(u_{j-1})\ts x^+_{j+1}(u_{j+1})\ldots x^+_{n}(u_{n}).
\eeq
\end{lem}

Regarding the identity \eqref{xm}, note that, in accordance with \eqref{notation788}, $x^-_0(u_0)$ is applied to the first and $x^+_{[n]}(u)$ on the next $n$ tensor factors of $\ndo\CC^N \ot (\ndo\CC^N)^{\ot n} \ot \hp$. Furthermore, in \eqref{xm}, as well as in the rest of the paper, we use the expansion convention  where the expressions of the form $(x_1+\ldots+x_n)^r$ with $r<0$ are expanded in the nonnegative powers of the variables $x_2,\ldots ,x_n$. Hence, for example, we have
$$
(u_0 - u_j)^r =\sum_{k\geqslant 0} \binom{r}{k} u_0^{r-k} (-u_j)^k
\in\CC[u_0^{-1}][[u_j]]\quad\text{for}\quad r<0.
$$

By employing \eqref{xm} one can prove
\beq\label{xmcom}
 x^-_1(u_1)  x^-_2(u_2)  =  x^-_2(u_2) x^-_1(u_1)  .
\eeq
Moreover, the equalities in \eqref{tr34} imply  
\beq\label{trxmcom}
x_{11}^{(r-1)}+x_{22}^{(r-1)}+\ldots + x_{NN}^{(r-1)}=0\quad\text{for all}\quad r=1,2,\ldots .
\eeq
We now regard $x^+(u)$ as an operator on $\hp$, where its action is given by the multiplication. By \eqref{xm} for $n=1$ we have
\beq\label{xmp}
 x^-_1 (u_1)  x^+_2(u_2) -x^+_2(u_2)  x^-_1 (u_1)= -S(u_1-u_2).
\eeq
Let us combine $x^+(u)$ and $x^-(u)$ into a single operator series
\beq\label{sglop}
x(u)=x^+(u)+x^-(u) \in  \ndo\CC^N \ot \om (\hp,\hp((u))_h).
\eeq
Since $P S(u) = S(u)P $, the identities \eqref{xmcom} and \eqref{xmp} imply
\beq\label{com12}
 x_1(u_1) x_2(u_2) -x_2(u_2) x_1(u_1) =-S (u_1-u_2)+S(u_2-u_1) ,
\eeq
while \eqref{tr23} and \eqref{trxmcom} imply
\beq\label{trcom12}
x_{11}^{(r)}+x_{22}^{(r)}+\ldots + x_{NN}^{(r)}=0\quad\text{for all}\quad r\in\ZZ.
\eeq
Note that equality \eqref{com12} can be also written as
\beq\label{com13}
x_1(u_1)x_2(u_2) + S (u_1-u_2) = x_2(u_2)x_1(u_1)+S(u_2-u_1).
\eeq
In accordance with our expansion convention,  the left hand side  of \eqref{com13} belongs to
$$
 (\ndo\CC^N)^{\ot 2}\ot \om(\hp,\hp ((u_1))((u_2))_h)   
$$
and the right hand side   to
$$
(\ndo\CC^N)^{\ot 2} \ot\om(\hp,\hp ((u_2))((u_1))_h),
$$
 so that the both sides  are elements of $(\ndo\CC^N)^{\ot 2}\ot\om(\hp,\hp ((u_1,u_2))_h)$. We write
\beq\label{normal2}
x_{[2]}(u)=x_{[2]}(u_1,u_2)= x_1(u_1)x_2(u_2) + S (u_1-u_2).
\eeq

Our next goal is to generalize \eqref{normal2} to an arbitrary number of factors. Let $n$ be a positive integer. For $n=1$ we set $x_{[1]}(u )=x(u )$. Suppose $n>1$. For any $k=0,\ldots,\floor{n/2} $ denote by $I_{k}^n$ the family of all sets of $k$   ordered pairs $(p,q)\in\left\{1,\ldots ,n\right\}^{\times 2}$ such that $p<q$ and such   that the coordinates   of all pairs   which belong to the same set are mutually distinct. For any $i\in I_{k}^n$  we denote by $i'$   the set of all integers in $\left\{1,\ldots ,n\right\}$ which do not appear in $i$.  

\begin{ex}
For $n=4$ and $k=0,1,2$ we have  $I_0^4 =\emptyset$,
\begin{gather*}
I_1^4 =\left\{ \left\{(1,2)\right\},\left\{(1,3)\right\}, \left\{(1,4)\right\},\left\{(2,3)\right\}, \left\{(2,4)\right\},\left\{(3,4)\right\}   \right\},\\
 I_2^4=\left\{\left\{(1,2),(3,4)  \right\}, \left\{(1,3),(2,4)  \right\}, \left\{(1,4),(2,3)  \right\}   \right\}.
\end{gather*}
Note that for all $i\in  I_2^4$ we have $i'=\emptyset$. As for the elements of $I_1^4$
we have, e.g., 
$$
\left\{(1,2)\right\}'=\left\{ 3,4 \right\},\quad \left\{(1,4)\right\}'=\left\{ 2,3 \right\},\quad\left\{(2,3)\right\}'=\left\{ 1,4 \right\}.
$$
\end{ex}

For   $u=(u_1,\ldots ,u_n)$ write $S_{i_p\ts j_p}=S_{i_p\ts j_p}(u_{i_p}-u_{j_p})$ for  all $p=1,\ldots ,k$. Define
\beq\label{normaln}
x_{[n]}(u) = 
\sum_{k=0}^{\floor{n/2}}
\sum_{ \substack{i=\left\{(i_1,j_1),\ldots ,(i_k,j_k)\right\}\in I^n_k\\l_1<\ldots< l_{n-2k} \in i' }}
S_{i_1\ts j_1} \ldots  S_{i_k\ts j_k} \ts x_{l_1}(u_{l_1})\ldots x_{l_{n-2k}}(u_{l_{n-2k}})
,
\eeq
where  the second sum  is $x_1(u_1)\ldots x_n(u_n)$ for $k=0$ and the indices denote the tensor factors of $(\ndo\CC^N)^{\ot n}\ot \hp$  on which the corresponding elements are applied. 

\begin{ex}
Clearly, \eqref{normaln} coincides with \eqref{normal2} for $n=2$, while for $n=3,4$ we get
\begin{align*}
x_{[3]}(u)=&\,\,x_1(u_1)\ts x_2(u_2)\ts x_3(u_3) + S_{12}\ts x_3(u_3)+S_{13}\ts x_2(u_2)+S_{23}\ts x_1(u_1),\\
x_{[4]}(u)=&\,\,x_1(u_1)\ts x_2(u_2)\ts x_3(u_3)\ts x_4(u_4) +S_{12}\ts x_3(u_3)\ts x_4(u_4) + S_{13}\ts x_2(u_2)\ts x_4(u_4) \\
&\,\,+ S_{14}\ts x_2(u_2)\ts x_3(u_3) + S_{23}\ts x_1(u_1)\ts x_4(u_4)
+ S_{24}\ts x_1(u_1)\ts x_3(u_3) \\
&\,\,+  S_{34}\ts x_1(u_1)\ts x_2(u_2) + S_{12}\ts S_{34} + S_{13}\ts S_{24}  + S_{14}\ts S_{23}. 
\end{align*}
\end{ex}

Obviously, \eqref{normaln} belongs to 
$(\ndo\CC^N)^{\ot n}\ot\om(\hp,\hp((u_1))\ldots ((u_n))_h)$. However, by using commutation relation \eqref{com13} one can verify the following stronger statement:

\begin{pro}\label{restricted}
For any  integer $n\geqslant 1$ we have
$$
x_{[n]}(u_1,\ldots ,u_n)\in (\ndo\CC^N)^{\ot n}\ot\om(\hp,\hp((u_1,\ldots ,u_n))_h).
$$
\end{pro}

Let $z$ be a single variable. Due to Proposition \ref{restricted}  we can define 
\beq\label{recall-}
x_{[n]}(z+u) = x_{[n]}(z+u_1,\ldots , z+u_n)\coloneqq
x_{[n]}(z_1,\ldots , z_n)\big|_{z_1=z+u_1,\ldots ,z_n=z+u_n}.\big.
\eeq
In addition, the given element satisfies
\beq\label{restricted2}
x_{[n]}(z+u) 
 \in 
(\ndo\CC^N)^{\ot n}\ot\om(\hp,\hp((z))_h [[u_1,\ldots ,u_n]]).
\eeq
Note that $x_{[n]}(z+u)$ differs from  
 $$
x_{[n]}(u+z) = x_{[n]}(u_1+z,\ldots , u_n+z)\coloneqq
x_{[n]}(z_1,\ldots , z_n)\big|_{z_1= u_1+z,\ldots ,z_n=u_n+z}\big. 
$$
(which is also well-defined by Proposition \ref{restricted})
as $x_{[n]}(u+z)$ should be expanded in the nonnegative powers of $z$.

\begin{rem}
As we   demonstrate later on, the coefficients of the matrix entries of  $x^+(u)$ and $x^-(u)$ can be regarded as   deformations of certain  creation and annihilation operators, respectively, while the operators $x_{[n]}(u)$ take place of the normal-ordered products. 
\end{rem}

\subsection{Braiding map}\label{subsec22}
From now on, the tensor products of $\CC[[h]]$-modules are understood as $h$-adically completed.
Define
\beq\label{te}
T(z)=S(z)-S(-z)\in\ndo\CC^N \ot\ndo\CC^N [z^{-1}]_h .
\eeq
Note that  
$
T(z_1 -z_2)$
is not equal to
$S (z_1 -z_2) - S(z_2-z_1)$, since $S(z_2-z_1)$ is to be expanded in nonnegative powers of $z_1$.
However, for any integer $n\geqslant 0 $ there exists an integer $r\geqslant 0$ such that
\beq\label{sloca}
(z_1 - z_2)^r \ts T(z_1 -z_2) =(z_1 - z_2)^r \left( S(z_1 -z_2) - S(z_2-z_1)\right)\mod h^n.
\eeq

Let $m,n\geqslant 1$ be integers.
For any $k=0,\ldots ,\min\left\{m,n\right\}$ let $I_{k}^{n,m}\subseteq I_{k}^{n+m} $ be the family of all sets of $k$   ordered pairs $(p,q)\in\left\{1,\ldots ,n\right\}\times \left\{n+1,\ldots ,n+m\right\}$   such that the coordinates of all pairs which
belong to the same set are mutually distinct. For   $i\in I_{k}^{n,m}$ we denote by $i'$ the set of all integers in $\left\{1,\ldots ,n+m\right\}$ which do not appear in $i$.
\begin{ex}
For $n=3$, $m=2$ and $k=0,1,2$ we have $I_0^{3,2}=\emptyset$,
\begin{align*}
I_1^{3,2}=&\left\{
\left\{(1,4)\right\}, \left\{(1,5)\right\}, \left\{(2,4)\right\}, \left\{(2,5)\right\}, \left\{(3,4)\right\}, \left\{(3,5)\right\}
\right\},\\
I_2^{3,2}=&\left\{
\left\{(1,4) , (2,5)\right\}, \left\{(1,4),(3,5)\right\}, \left\{(1,5),(2,4)\right\}, \left\{(1,5),(3,4)\right\},\right.\\
&\,\, \left.\left\{(2,4),(3,5)\right\}, \left\{(2,5),(3,4)\right\} 
\right\}.
\end{align*}
Also, for example, we have 
$$
\left\{(1,4)\right\}'=\left\{2,3,5\right\},\quad
\left\{(3,5)\right\}'=\left\{1,2,4\right\},\quad
\left\{(2,4),(3,5)\right\}'=\left\{1\right\}.
$$
\end{ex}

Let $u=(u_1,\ldots ,u_n)$, $v=(v_1,\ldots ,v_m)$ be families of variables. Consider the expression
\beq\label{tensors2}
x_{[n]}^{+13} (u)x_{[m]}^{+24} (v)
=x_{1\ts n+m+1}^+ (u_1)\ldots x_{n\ts n+m+1}^+ (u_n)\ts
x_{n+1\ts n+m+2}^+ (v_1)\ldots x_{n+m\ts n+m+2}^+ (v_m)
\eeq
with coefficients in
\beq\label{tensors}
\smalloverbrace{(\ndo\CC^N)^{\ot n} }^{1}\ot
\smalloverbrace{(\ndo\CC^N)^{\ot m}  }^{2}\ot
 \smalloverbrace{\hp}^{3} \ot
 \smalloverbrace{\hp}^{4}
\eeq
and superscripts $1,2,3,4$ indicating the tensor factors as in \eqref{tensors}. If $n=0$ or $m=0$ we define the corresponding empty product in \eqref{tensors2} to be the unit $\vac\in\hp$.

\begin{lem}\label{blemma}
There exists a unique $\CC[[h]]$-module map $$\Sc(z)\colon \hp\ot\hp\to \hp\ot\hp[z^{-1}]_h$$ 
such that for any integers $m,n\geqslant 0$ and the variables $u=(u_1,\ldots, u_n)$ and $v=(v_1,\ldots ,v_m)$
\begin{align}
\Sc(z)\ts x_{[n]}^{+13} (u)x_{[m]}^{+24} (v)  =
\sum_{k=0}^{\min\left\{m,n\right\}}
\sum_{ \substack{i=\left\{(i_1,j_1),\ldots ,(i_k,j_k)\right\}\in I_k^{n,m}\\l_1<\ldots< l_{n+m-2k} \in i' }}
T_{i_1\ts j_1} \ldots  T_{i_k\ts j_k} \ts x^+_{l_1} \ldots x^+_{l_{n+m-2k}} ,\label{es}
\end{align}
where    $T_{i_p\ts j_p} $ denotes  $T_{i_p\ts j_p} (z+u_{i_p}-v_{j_p -n})$ and
$$
x^+_{l_p}=\begin{cases}
x^+_{l_p\ts n+m+1}(u_{l_p})&\text{for }l_p=1,\ldots ,n,\\
x^+_{l_p\ts n+m+2}(v_{l_p -n})&\text{for }l_p=n+1,\ldots ,n+m .
\end{cases}
$$
The given map   is of the form $\Sc=1+O(h)$ and satisfies the {\em Yang--Baxter equation}, 
\beq\label{s2}
\mathcal{S}_{12}(z_1)\ts\mathcal{S}_{13}(z_1+z_2)\ts\mathcal{S}_{23}(z_2)
=\mathcal{S}_{23}(z_2)\ts\mathcal{S}_{13}(z_1+z_2)\ts\mathcal{S}_{12}(z_1) 
\eeq
and the {\em unitarity condition},  
\beq\label{s3}
\mathcal{S}_{21}(z)=\mathcal{S}^{-1}(-z).
\eeq
\end{lem}

\begin{prf}
It is clear that \eqref{es} uniquely determines a well-defined $\CC[[h]]$-module map on $\hp\ot\hp$. Furthermore, its image belongs to $\hp\ot\hp[z^{-1}]_h$ as $T(z)\in\ndo\CC^N\ot\ndo\CC^N [z^{-1}]_h$. In fact, \eqref{esform} and \eqref{te} imply that $T(z)$ belongs to $h\ndo\CC^N\ot\ndo\CC^N [z^{-1}]_h$. Hence,  as the only summand on the right-hand side of \eqref{es} which does not contain a copy of $T$ is $x_{[n]}^{+13} (u)x_{[m]}^{+24} (v) $, the given map is of the form $\Sc=1+O(h)$ . 

Regarding the Yang--Baxter equation \eqref{s2}, it is sufficient to show that the operators $\mathcal{S}_{12}(z_1)$, $\mathcal{S}_{13}(z_1+z_2)$ and $\mathcal{S}_{23}(z_2)$ are mutually commutative. However, this follows directly from the definition \eqref{es}. Indeed, choose
any   $A,B\in\left\{\mathcal{S}_{12}(z_1) ,  \mathcal{S}_{13}(z_1+z_2), \mathcal{S}_{23}(z_2)\right\}$ such that $A\neq B$ and then consider the actions of $AB$ and $BA$
on the expression
$$
x\coloneqq x_{[n]}^{+14} (u)x_{[m]}^{+25} (v)x_{[k]}^{+36} (w)
=
x_{[n]}^{+14} (u_1,\ldots,u_n)x_{[m]}^{+25} (v_1,\ldots ,v_m)x_{[k]}^{+36} (w_1,\ldots ,w_k) 
$$
with coefficients in
\beq\label{bhj}
(\ndo\CC^N)^{\ot n }\ot(\ndo\CC^N)^{\ot m}\ot(\ndo\CC^N)^{\ot k}\ot \hp \ot \hp \ot \hp .
\eeq
By \eqref{es}, both of the  actions produce a sum of certain elements of the form
\beq\label{fctrs}
T_{i_1\ts j_1}\ldots T_{i_r\ts j_r}\ts x^+_{a_1\ts b_1}\ldots x^+_{a_{n+m+k-2r}\ts b_{n+m+k-2r}}
\eeq
with coefficients in \eqref{bhj} such that all factors in \eqref{fctrs}  mutually commute. Finally,  note that every such term of the form \eqref{fctrs} appears as a summand in $ABx$ if and only if it appears as a summand in $BAx$, which implies $AB=BA$.

Let us prove that $\Sc(z)$ possesses the unitarity property \eqref{s3}.
Note that $T(z)$ is an odd function, i.e. we have $T(z)=-T(-z)$; recall \eqref{te}. Moreover, as $R_{12}(z)=R_{21}(z)$ we have  $T_{12}(z)=T_{21}(z)$ and, consequently, $\Sc_{12}(z)=\Sc_{21}(z)$. Therefore, by \eqref{es} we have
\begin{align}
\Sc_{21}(-z)\ts x_{[n]}^{+13} (u)x_{[m]}^{+24} (v)  =
\sum_{k=0}^{\min\left\{m,n\right\}}
\sum_{ \substack{i=\left\{(i_1,j_1),\ldots ,(i_k,j_k)\right\}\in I^{n,m}_k\\l_1<\ldots< l_{n+m-2k} \in i' }}(-1)^k
T_{i_1\ts j_1} \ldots  T_{i_k\ts j_k} \ts x^+_{l_1} \ldots x^+_{l_{n+m-2k}}.\label{es2}
\end{align}
By combining the identities \eqref{es} and \eqref{es2} one easily checks that all terms of the form 
$T_{i_1\ts j_1} \ldots  T_{i_k\ts j_k} \ts x^+_{l_1} \ldots x^+_{l_{n+m-2k}}$ with $k>0$ in  
\beq\label{esovi}
\Sc_{21}(-z)\ts \Sc(z)\ts x_{[n]}^{+13} (u)x_{[m]}^{+24} (v) \fand
\Sc(z)\ts\Sc_{21}(-z)\ts  x_{[n]}^{+13} (u)x_{[m]}^{+24} (v)
\eeq
cancel, so that both expressions in \eqref{esovi} are equal to
$ x_{[n]}^{+13} (u)x_{[m]}^{+24} (v) $,
as required.
\end{prf}

\subsection{Quantum vertex algebra \texorpdfstring{$\hp$}{HN+}}\label{subsec23}

Let us    recall the Etingof--Kazhdan  definition of quantum vertex algebra; see \cite[Subsect. 1.4]{EK}.
\begin{defn}\label{qvoa}
A {\em quantum vertex algebra} is a quadruple $(V,Y,\vac,\Sc)$ such that
\begin{enumerate}[(1)]
\item  $V$ is a topologically free $\mathbb{C}[[h]]$-module.
\item $Y$ is a $\mathbb{C}[[h]]$-module map (the {\em vertex operator map})
\begin{align*}
Y \colon V\ot V&\to V((z))_h\\
u\ot v&\mapsto Y(z)(u\ot v)=Y(u,z)v=\sum_{r\in\mathbb{Z}} u_r v \ts z^{-r-1}
\end{align*}
which satisfies the {\em weak associativity}:
for any $u,v,w\in V$ and $n\in\mathbb{Z}_{\geqslant 0}$
there exists $s\in\mathbb{Z}_{\geqslant 0}$
such that
\begin{equation}\label{associativity}
(z_0 +z_2)^s\ts Y(u,z_0 +z_2)Y(v,z_2)\ts w - (z_0 +z_2)^s\ts Y\big(Y(u,z_0)v,z_2\big)\ts w
\in h^n V[[z_0^{\pm 1},z_2^{\pm 1}]].
\end{equation}
\item $\vac$ is a distinct element of $V$ (the {\em vacuum vector}) such that
\beq\label{v1}
Y(\vac ,z)v=v,
\quad 
Y(v,z)\ts\vac  \in V[[z]]\fand \lim_{z\to 0} Y(v,z)\ts\vac =v\qquad\text{for all }v\in V.
\eeq
\item $\Sc=\Sc(z)$ is a $\mathbb{C}[[h]]$-module map
$ V\otimes V\to V\otimes V\otimes\mathbb{C}((z))[[h]]$ (the {\em braiding}) of the form $\Sc=1+O(h)$ which satisfies 
the Yang--Baxter equation \eqref{s2}, the unitarity condition \eqref{s3},
the {\em shift condition}
\begin{align}
&[D\otimes 1, \mathcal{S}(z)]=-\frac{d}{dz}\mathcal{S}(z),\quad\text{where}\quad Dv\coloneqq v_{-2}\vac\quad\text{for all }v\in V,\label{s1}
\end{align}
 the $\mathcal{S}$-{\em locality}:
for any $u,v\in V$ and $n\in\mathbb{Z}_{\geqslant 0}$ there exists
$r\in\mathbb{Z}_{\geqslant 0}$ such that  
\begin{align}
&(z_1-z_2)^{r}\ts Y(z_1)\big(1\otimes Y(z_2)\big)\big(\mathcal{S}(z_1 -z_2)(u\otimes v)\otimes w\big)
\nonumber\\
&\quad-(z_1-z_2)^{r}\ts Y(z_2)\big(1\otimes Y(z_1)\big)(v\otimes u\otimes w)
\in h^n V[[z_1^{\pm 1},z_2^{\pm 1}]]\quad\text{for all }w\in V\label{locality}
\end{align}
and the {\em hexagon identity}:
\beq\label{hex1}
\Sc(z_1)\left(Y(z_2) \ot 1\right) = \left(Y(z_2) \ot 1\right) \Sc_{23}(z_1)\Sc_{13}(z_1+z_2).
\eeq
\end{enumerate}
\end{defn}

\begin{rem}\label{hexrem}
We should mention that Definition \ref{qvoa} slightly differs from the  original   in \cite{EK}. We included both   weak associativity \eqref{associativity} and   hexagon identity \eqref{hex1} in our definition in order  to emphasize the importance of both of these properties for quantum vertex algebra theory. However, any one of them can be omitted from the definition and then proved using the remaining axioms; see \cite[Prop. 1.4]{EK} and \cite[Prop. 3.14]{DGK}.
\end{rem}
 
We now use the data from the previous subsections, in particular, Lemmas \ref{lemma21} and   \ref{blemma}, to define a quantum vertex algebra structure over $\hp$.

\begin{thm}\label{mainthm}
For any $c\in\CC$ there exists a unique structure of quantum vertex algebra on $\hp$ such that the vacuum vector is the unit $\vac\in\hp$, the braiding $\Sc(z)$ is given by \eqref{es}    and the vertex operator map is given by
\beq\label{ymap}
Y(x_{[n]}^+ (u_1,\ldots ,u_n),z)
=
x_{[n]}(z+u_1,\ldots , z+u_n).
\eeq
\end{thm}

The theorem is proved by directly verifying the constraints imposed by Definition \ref{qvoa}. We dedicate the next subsection to its proof.

\subsection{A   proof of Theorem \ref{mainthm}}\label{appendix}\label{subsec24}

It is clear that the $h$-adically completed polynomial algebra $\hp$ is topologically free. Moreover, the matrix entries of the coefficients of all $x_{[n]}^+ (u_1,\ldots ,u_n)$ along with $\vac$ span an $h$-adically dense $\CC[[h]]$-submodule of $\hp$ so that \eqref{ymap} uniquely determines the vertex operator map. However, we have to prove that the vertex operator map is well-defined by \eqref{ymap}. It is sufficient to check that $Y(z)$ maps the ideal of defining relations
$
\left[x_{1}^+(u),x_2^{+}(v)\right]=0 
$
and \eqref{tr23} for $\hp$ to itself. Regarding the first family of relations, for any integers $n>j>0$ consider the expression
\beq\label{img3}
x_{[n]}^{+(j)}(u)\coloneqq x_1^+(u_1)\ldots x_{j-1}^+(u_{j-1})\ts x_{j+1}^+(u_{j+1})\ts x_j^+(u_j)\ts x_{j+2}^+(u_{j+2})
\ldots x_n^+(u_n).
\eeq
 We will show that its image under $Y(z)$ coincides with the right-hand side of \eqref{ymap}, which implies the desired conclusion.  First, we write \eqref{normaln} as
\beq\label{wdf1}
x_{[n]}(u)=x_{[n]}(u)_{j,j+1} + x_{[n]}(u)^{j,j+1} + x_{[n]}(u)_{0},\qquad\text{where}
\eeq
\begin{itemize}
\item $x_{[n]}(u)_{j,j+1}$ denotes the sum of all $S_{i_1\ts j_1} \ldots  S_{i_k\ts j_k} \ts x_{l_1}(u_{l_1})\ldots x_{l_{n-2k}}(u_{l_{n-2k}})$ in \eqref{normaln} which contain $x_j(u_j)x_{j+1}(u_{j+1})$;
\item $x_{[n]}(u)^{j,j+1}$ denotes the sum of all $S_{i_1\ts j_1} \ldots  S_{i_k\ts j_k} \ts x_{l_1}(u_{l_1})\ldots x_{l_{n-2k}}(u_{l_{n-2k}})$ in \eqref{normaln} which contain $S_{j\ts j+1}(u_j - u_{j+1})$;
\item  $x_{[n]}(u)_{0}= x_{[n]}(u)-x_{[n]}(u)_{j,j+1} - x_{[n]}(u)^{j,j+1}$ are the remaining terms.
\end{itemize}
Next, we  compare $x_{[n]}(u)$ with
$$
x_{[n]}^{(j)}(u)\coloneqq P_{j\ts j+1} \ts x_{[n]}(u_1,\ldots ,u_{j-1},u_{j+1},u_j,u_{j+2},\ldots ,u_n)\ts P_{j\ts j+1},
$$ 
where $P_{j\ts j+1}$ denotes the permutation operator from \eqref{perm8} applied on the tensor factors $j$ and $j+1$.
As with \eqref{wdf1}, we write $x_{[n]}^{(j)}(u)$ as
\beq\label{wdf2}
x_{[n]}^{(j)}(u)=x_{[n]}^{(j)}(u)_{j+1,j} + x_{[n]}^{(j)}(u)^{j+1,j} + x_{[n]}^{(j)}(u)_{0},\qquad\text{where}
\eeq
\begin{itemize}
\item $x_{[n]}^{(j)}(u)_{j+1,j}$ denotes the sum of all $S_{i_1\ts j_1} \ldots  S_{i_k\ts j_k} \ts x_{l_1}(u_{l_1})\ldots x_{l_{n-2k}}(u_{l_{n-2k}})$ in $x_{[n]}^{(j)}(u)$ which contain $x_{j+1}(u_{j+1})x_j(u_j)$;
\item $x_{[n]}^{(j)}(u)^{j+1,j}$ denotes the sum of all $S_{i_1\ts j_1} \ldots  S_{i_k\ts j_k} \ts x_{l_1}(u_{l_1})\ldots x_{l_{n-2k}}(u_{l_{n-2k}})$ in $x_{[n]}^{(j)}(u)$ which contain $S_{j+1\ts j}(u_{j+1} - u_j)$;
\item $x_{[n]}^{(j)}(u)_{0}= x_{[n]}^{(j)}(u)-x_{[n]}^{(j)}(u)_{j+1,j} - x_{[n]}^{(j)}(u)^{j+1,j}$ are the remaining terms.
\end{itemize}
Clearly, $x_{[n]}(u)_{0}= x_{[n]}^{(j)}(u)_{0}$. Furthermore, as $S_{12}(z)=S_{21}(z)$, the identity \eqref{com13} implies 
$$
x_{[n]}(u)_{j,j+1} + x_{[n]}(u)^{j,j+1} = x_{[n]}^{(j)}(u)_{ j+1,j}+x_{[n]}^{(j)}(u)^{ j+1,j}.
$$
By combining these two observations with \eqref{wdf1} and \eqref{wdf2} we find that
\beq\label{wdf3}
x_{[n]}(u)=x_{[n]}^{(j)}(u).
\eeq
Finally, by replacing the variables $(u_1,\ldots ,u_n)$ with $(z+u_1,\ldots ,z+u_n)$ in \eqref{wdf3} we get
$$
Y(x_{[n]}^+(u), z) =Y(x_{[n]}^{+(j)}(u),z),
$$
as required. As for the remaining family of relations \eqref{tr23},  this is a direct consequence of \eqref{tr34} and \eqref{trcom12}.
Indeed, due to the aforementioned equalities, for any $i=1,\ldots ,n$ by applying the partial trace $\tr_i$ to the $i$-th tensor factor of any summand in \eqref{normaln} produces zero.
 Thus, we conclude that the vertex operator map is well-defined by \eqref{ymap}.

Regarding  the vacuum vector axioms \eqref{v1}, $Y(\vac ,z)v=v$ clearly holds for all $v\in \hp$ while the remaining properties follow from the identity
\beq\label{acvac}
x_{[n]}(u_1,\ldots ,u_n)\vac= x_{[n]}^+ (u_1,\ldots ,u_n), 
\eeq
which can be directly verified using Lemma \ref{lemma21} and commutation  relations \eqref{xmp}. Finally,   the image of the vertex operator map belongs to $\hp\ot\hp ((z))_h$   by \eqref{restricted2}. 

As for the braiding, let us verify the shift condition \eqref{s1}. The remaining requirements on $\Sc(z)$, which are imposed by Definition \ref{qvoa}, hold by Lemma \ref{blemma}. First, note that by applying \eqref{ymap} on the vacuum vector $\vac$ and then taking the coefficient with respect to $z$ we get, due to \eqref{acvac},
$$
D\ts x_{[n]}^+(u_1,\ldots ,u_n) = \textstyle\left(\sum_{r=1}^n \frac{\partial}{\partial u_r}\right) x_{[n]}^+(u_1,\ldots ,u_n).
$$
Furthermore, by \eqref{v1} we have $D\ts\vac =0$. Using these observations along with \eqref{es} we compute the action of $(D\ot 1)  \Sc(z)$ and $\Sc(z) (D\ot 1)  $ on
$$
x\coloneqq x_{[n]}^{+13} (u)x_{[m]}^{+24} (v) 
=
x_{[n]}^{+13} (u_1,\ldots,u_n)x_{[m]}^{+24} (v_1,\ldots ,v_m) 
$$
as follows. We get
\begin{align}
&(D\ot 1)\Sc(z)\ts x = (D\ot 1)
\sum_{k=0}^{\min\left\{m,n\right\}}
\sum_{ \substack{i=\left\{(i_1,j_1),\ldots ,(i_k,j_k)\right\}\in I_k^{n,m}\\l_1<\ldots< l_{n+m-2k} \in i' }}
T_{i_1\ts j_1} \ldots  T_{i_k\ts j_k} \ts x^+_{l_1} \ldots x^+_{l_{n+m-2k}}\non \\
&\qquad=  
\sum_{k=0}^{\min\left\{m,n\right\}}
\sum_{ \substack{i=\left\{(i_1,j_1),\ldots ,(i_k,j_k)\right\}\in I_k^{n,m}\\l_1<\ldots< l_{n+m-2k} \in i' }}
T_{i_1\ts j_1} \ldots  T_{i_k\ts j_k} \ts \textstyle\left(\sum_{r=1}^n \frac{\partial}{\partial u_r}\right)\left( x^+_{l_1} \ldots x^+_{l_{n+m-2k}}\right),\label{hq1}\\
&\Sc(z) (D\ot 1)\ts x=\Sc(z)\textstyle\left(\sum_{r=1}^n \frac{\partial}{\partial u_r}\right) x_{[n]}^{+13} (u_1,\ldots,u_n)  x_{[m]}^{+24} (v_1,\ldots ,v_m)\non \\
&\qquad=
\textstyle\left(\sum_{r=1}^n \frac{\partial}{\partial u_r}\right) \displaystyle
\sum_{k=0}^{\min\left\{m,n\right\}}
\sum_{ \substack{i=\left\{(i_1,j_1),\ldots ,(i_k,j_k)\right\}\in I_k^{n,m}\\l_1<\ldots< l_{n+m-2k} \in i' }}
T_{i_1\ts j_1} \ldots  T_{i_k\ts j_k} \ts x^+_{l_1} \ldots x^+_{l_{n+m-2k}},\label{hq2}
\end{align}
where all $T_{i_p\ts j_p} $ and $x^+_{l_p}$ are defined as in the statement of Lemma \ref{blemma}.
Finally, by comparing the expressions \eqref{hq1} and \eqref{hq2} and using
the identity
$$
\textstyle\left(\sum_{r=1}^n \frac{\partial}{\partial u_r}\right)
T_{i_1\ts j_1} \ldots  T_{i_k\ts j_k}
=\frac{\partial}{\partial z} T_{i_1\ts j_1} \ldots  T_{i_k\ts j_k}
,
$$
which follows from
$$
\textstyle\frac{\partial}{\partial u_{i_p}} T_{i_p\ts j_p} (z+u_{i_p}-v_{j_p -n})
=\frac{\partial}{\partial z} T_{i_p\ts j_p} (z+u_{i_p}-v_{j_p -n}),
$$
we conclude that the difference of \eqref{hq1} and \eqref{hq2} is equal to $-\frac{\partial}{\partial z}\Sc(z)\ts x $, as required.
Hence the shift condition holds.

To finish the proof it remains to check that the vertex operator map satisfies the
 weak associativity \eqref{associativity} and $\Sc$-locality \eqref{locality}; recall Remark \ref{hexrem}.  These    properties are verified in Lemmas \ref{assocl} and \ref{locall}  below. Alternatively, one could  prove   the hexagon identity \eqref{hex1} instead of the weak associativity.  However, we find the  direct  proof of the hexagon identity to be more technical, so   we prove the latter property instead.

Let $u=(u_1,\ldots ,u_n)$ be a family of variables and $z$ a single variable. In order to simplify the notation, we denote the families $(z+u_1,\ldots ,z+u_n)$ and $(u_1+z,\ldots ,u_n+z)$ by $z+u$ and $u+z$ respectively. For example, \eqref{ymap} can be written briefly as
$$Y(x_{[n]}^+ (u ),z)
=
x_{[n]}(z+u).$$

\begin{lem}\label{assocl}
Vertex operator map  \eqref{ymap} satisfies the weak associativity   \eqref{associativity}.
\end{lem}

\begin{prf}
Choose any element $y\in\hp$ and integers $t,r_1,\ldots ,r_n,s_1,\ldots ,s_m\in\ZZ_{\geqslant 0}$. 
We start by considering the image of
\beq\label{asso1}
x_{[n]}^{+13}(u)\ts x_{[m]}^{+24}(v)=x_{1\ts n+m+1}^+(u_1)\ldots x_{n\ts n+m+1}^+(u_n)\ts x_{n+1\ts n+m+2}^+(v_1)\ldots x_{n+m\ts n+m+2}^+(v_m)
\eeq
under the second summand in \eqref{associativity}. 
By applying $Y(z_2)\left(Y(z_0)\ot 1\right)$ to \eqref{asso1} we get
\beq\label{asso3}
Y (x_{[n]}^{13}(z_0+u)  x_{[m]}^{+23}(v), z_2 ).
\eeq
By \eqref{xmp} the argument of the vertex operator map can be expressed as
\begin{align}\label{asso2}
x_{[n]}^{13}(z_0+u)  x_{[m]}^{+23}(v)=
\sum_{k=0}^{\min\left\{m,n\right\}}
\sum_{ \substack{i=\left\{(i_1,j_1),\ldots ,(i_k,j_k)\right\}\in I^{n,m}_k\\l_1<\ldots< l_{n+m-2k} \in i' }}(-1)^k
S_{i_1\ts j_1} \ldots  S_{i_k\ts j_k} \ts x^+_{l_1} \ldots x^+_{l_{n+m-2k}} 
\end{align}
for $S_{i_p\ts j_p}=S_{i_p\ts j_p}(z_0+u_{i_p}-v_{j_p -n})=S_{i_p\ts j_p}(w_{i_p}-w_{j_p})$ and $x^+_{l_p}=x^+_{l_p\ts n+m+1}(w_{l_p})$, where
$$
w_{l_p}=\begin{cases}
 z_0 +u_{l_p}&\text{for }l_p=1,\ldots ,n,\\
 v_{l_p -n}&\text{for }l_p=n+1,\ldots ,n+m.
\end{cases}
$$
Therefore, by applying $Y(z_2)$ to \eqref{asso2} and using \eqref{ymap} we conclude that \eqref{asso3} equals
\beq\label{asso6}
\sum_{k=0}^{\min\left\{m,n\right\}}
\sum_{ \substack{i=\left\{(i_1,j_1),\ldots ,(i_k,j_k)\right\}\in I_k^{n,m}\\l_1<\ldots< l_{n+m-2k} \in i' }}(-1)^k
S_{i_1\ts j_1} \ldots  S_{i_k\ts j_k} \ts   x_{[n+m-2k]} (z_2+w_{l_1},\ldots  ,z_2+w_{l_{n+m-2k}}),
\eeq
where the expression
$
x_{[n+m-2k]} (z_2+w_{l_1},\ldots  ,z_2+w_{l_{n+m-2k}})
$
is applied on the tensor factors $l_1,\ldots ,l_{n+m-2k}$ and $n+m+1$ of
\beq\label{asso5}
(\ndo\CC^N)^{\ot n}\ot (\ndo\CC^N)^{\ot m} \ot \hp.
\eeq
Define
$$
x_{l_p}=\begin{cases}
 z_0 +z_2+u_{l_p}&\text{for }l_p=1,\ldots ,n,\\
 z_2+v_{l_p -n}&\text{for }l_p=n+1,\ldots ,n+m.
\end{cases}
$$
Clearly,  there are only finitely many summands in \eqref{asso6}. Hence by Proposition \ref{restricted} there exists a nonnegative integer $s$ such that the coefficients of  
\beq\label{varsasso}
u_1^{r'_1}\ldots u_n^{r'_n} v_1^{s'_1}\ldots v_m^{s'_m}\quad\text{with}\quad r'_i=0,\ldots,r_i,\fand s'_j=0,\ldots ,s_j,
\eeq
where $i=1,\ldots ,n$ and $j=1,\ldots ,m$,
in  
\begin{align}
&(z_0 + z_2)^s S_{i_1\ts j_1} \ldots  S_{i_k\ts j_k} \ts   x_{[n+m-2k]} (z_2+w_{l_1},\ldots  ,z_2+w_{l_{n+m-2k}})\ts y\fand\label{ztr1}\\
&(z_0 + z_2)^s S_{i_1\ts j_1} \ldots  S_{i_k\ts j_k} \ts   x_{[n+m-2k]} (x_{l_1},\ldots  ,x_{l_{n+m-2k}})\ts y\label{ztr2}
\end{align}
coincide modulo $h^t$. Note that $x_{l_p}=z_2+w_{l_p}$. However, the variable $z_0$ comes   first   from the left in $x_{l_p}$ for $l_p=1,\ldots,n$, thus indicating that different expansions are applied   in \eqref{ztr1} and \eqref{ztr2}. Hence the expressions in \eqref{ztr1} and \eqref{ztr2} do not need to be   equal.

For any choice of indices  $l_1,\ldots ,l_{n+m-2k}$ as in \eqref{asso6} introduce the function  
\begin{align*}
\sigma\coloneqq \sigma_{l_1,\ldots ,l_{n+m-2k}}\colon\left\{1,\ldots,n+m-2k\right\}&\to \left\{l_1,\ldots,l_{n+m-2k}\right\}\\
p&\mapsto l_p .
\end{align*}
Henceforth all the given expressions are regarded modulo $u_{1}^{r_1 +1}\ldots u_{n}^{r_n +1} v_1^{s_1 +1}\ldots v_{m}^{s_m +1}h^t$, i.e. we work with the coefficients of \eqref{varsasso} modulo $h^t$.
Consider the action of the product of $(z_0 + z_2)^s$ and \eqref{asso6} on the element $y$.
By rewriting all terms $
x_{[n+m-2k]} (z_2+w_{l_1},\ldots  ,z_2+w_{l_{n+m-2k}})
$ of \eqref{asso6} using   formula \eqref{normaln} we get
\begin{align}
&(z_0 + z_2)^s \sum_{k=0}^{\min\left\{m,n\right\}}
\sum_{ \substack{i=\left\{(i_1,j_1),\ldots ,(i_k,j_k)\right\}\in I_k^{n,m}\\l_1<\ldots< l_{n+m-2k} \in i' }}(-1)^k
S_{i_1\ts j_1} \ldots  S_{i_k\ts j_k}  \non\\
&\times
\sum_{r=0}^{\floor{(n+m-2k)/2}}
\sum_{ \substack{a=\left\{(a_1,b_1),\ldots ,(a_{r},b_{r})\right\}\in \sigma(I^{n+m-2k}_{r})\\c_1<\ldots< c_{n+m-2k-2r}\ts \in\ts \sigma(a') }}
S_{a_1\ts b_1} \ldots  S_{a_{r}\ts b_{r}} 
\ts x_{c_1} \ldots x_{c_{n+m-2k-2r}} \ts y, \label{asso7}
\end{align}
where $\sigma= \sigma_{l_1,\ldots ,l_{n+m-2k}}$, 
$ 
S_{a_p\ts b_p}=S_{a_p\ts b_p}(w_{a_p}-w_{b_p})$ and
$x_{c_p}=x_{c_p\ts n+m+1}(z_2 +w_{c_p})$.
All summands in \eqref{asso7} are of the form
\beq\label{asso8}
(-1)^k
S_{i_1\ts j_1} \ldots  S_{i_k\ts j_k}\ts S_{a_1\ts b_1} \ldots  S_{a_{r}\ts b_{r}} 
\ts x_{c_1} \ldots x_{c_{n+m-2k-2r}}\ts y.
\eeq
However, due to the sign $(-1)^k$, all summands \eqref{asso8} in \eqref{asso7} which contain at least one copy of $S_{p\ts q}$ with $p\in\left\{1,\ldots ,n\right\}$ and $q\in\left\{n+1,\ldots ,n+m\right\}$ cancel, so that \eqref{asso7}   equals
\begin{align}
(z_0 + z_2)^s
\sum_{r=0}^{\floor{(n+m)/2}}
\sum_{ \substack{a=\left\{(a_1,b_1),\ldots ,(a_{r},b_{r})\right\}\in  I^{n+m}_{r}\setminus I^{n,m}_{r} \\c_1<\ldots< c_{n+m-2r}\ts \in\ts  a'  }}
S_{a_1\ts b_1} \ldots  S_{a_{r}\ts b_{r}} 
\ts x_{c_1} \ldots x_{c_{n+m-2r}}\ts y.\label{asso9}
\end{align}
Finally, we conclude by \eqref{normaln} that \eqref{asso9}  coincides with
\beq\label{assoa}
 (z_0 + z_2)^s\ts  x_{[n]}^{13}(z_2+z_0+u)\ts x_{[m]}^{23}(z_2+v)\ts y.
\eeq

Now consider the image of \eqref{asso1} under the first summand in \eqref{associativity}.
By applying $ (z_0 + z_2)^s\ts Y(z_0+z_2)(1\ot  Y(z_2)) $ to \eqref{asso1} and using \eqref{ymap} we get
\beq\label{assob}
 (z_0 + z_2)^s\ts x_{[n]}^{13}(z_0+z_2+u)\ts x_{[m]}^{23}(z_2+v)\ts y.
\eeq
By the choice of the integer $s$, the coefficients of the variables \eqref{varsasso} in \eqref{assoa} and \eqref{assob} coincide modulo $h^t$, so the weak associativity follows.
\end{prf}

\begin{lem}\label{locall}
Vertex operator map  \eqref{ymap} possesses the $\Sc$-locality property \eqref{locality}.
\end{lem}

\begin{prf}
As with the proof of Lemma \ref{assocl}, we verify the $\Sc$-locality directly, by comparing the images of \eqref{asso1} under the first and the second summand in \eqref{locality}. 
Let $t,r_1,\ldots ,r_n,s_1,\ldots ,s_m \geqslant 0 $ be arbitrary integers. 
Henceforth all the   expressions are regarded modulo 
\beq\label{modulo97}
u_{1}^{r_1 +1}\ldots u_{n}^{r_n +1} v_1^{s_1 +1}\ldots v_{m}^{s_m +1}h^t,
\eeq 
i.e. we consider only the coefficients of the monomials \eqref{varsasso} modulo $h^t$.
Applying the first summand in \eqref{locality}, $Y(z_1)\big(1\otimes Y(z_2)\big) \mathcal{S}(z_1 -z_2) $ to \eqref{asso1}
 and   using \eqref{es} we get
\begin{align}
Y(z_1)\big(1\otimes Y(z_2)\big)
\sum_{k=0}^{\min\left\{m,n\right\}}
\sum_{ \substack{i=\left\{(i_1,j_1),\ldots ,(i_k,j_k)\right\}\in I_k^{n,m}\\l_1<\ldots< l_{n+m-2k} \in i' }}
T_{i_1\ts j_1} \ldots  T_{i_k\ts j_k} \ts x^+_{l_1} \ldots x^+_{l_{n+m-2k}},\label{sloc1}
\end{align}
where   
\begin{align}
&T_{i_p\ts j_p}   
= S_{i_p\ts j_p} (z_1 - z_2+u_{i_p}-v_{j_p -n})
-S_{i_p\ts j_p} (-z_1 + z_2-u_{i_p}+v_{j_p -n}),\label{sloc8}\\
&x^+_{l_p}=\begin{cases}
x^+_{l_p\ts n+m+1}(u_{l_p})&\text{for }l_p=1,\ldots ,n,\\
x^+_{l_p\ts n+m+2}(v_{l_p -n})&\text{for }l_p=n+1,\ldots ,n+m .
\end{cases}\non
\end{align}
Next, using \eqref{ymap}, we rewrite \eqref{sloc1} as
\beq\label{sloc2}
\sum_{k=0}^{\min\left\{m,n\right\}}
\sum_{ \substack{i=\left\{(i_1,j_1),\ldots ,(i_k,j_k)\right\}\in I_k^{n,m}\\l_1<\ldots< l_{n+m-2k} \in i' }}
T_{i_1\ts j_1} \ldots  T_{i_k\ts j_k} \ts x_{[n-k]}^{13}\ts x_{[m-k]}^{23},
\eeq
where the terms
\begin{align}
&x_{[n-k]}^{13}=x_{[n-k]}(z_1+u_{l_1},\ldots ,z_1+u_{l_{n-k}})\Fand\label{sloc5}
\\
 &x_{[m-k]}^{23}=x_{[m-k]} (z_2+v_{l_{n-k+1}-n},\ldots , z_2+v_{l_{n+m-2k}-n})\label{sloc6}
\end{align}
are applied on the tensor factors
\beq\label{sloc7}
l_1,\ldots, l_{n-k},\, n+m+1\Fand
l_{n-k+1},\ldots, l_{n+m-2k},\, n+m+1
\eeq
of \eqref{asso5}, respectively.

For every $T_{i_p\ts j_p} $, as given by \eqref{sloc8}, introduce the element 
\beq\label{sloc9}
U_{i_p\ts j_p}=S_{i_p\ts j_p}(z_1 -z_2 +u_{i_p}-v_{j_p -n}) -S_{i_p\ts j_p}(z_2 -z_1+v_{j_p -n} - u_{i_p}).
\eeq
It is clear that   $T_{i_p\ts j_p} $  and $U_{i_p\ts j_p}$ do not coincide due to different expansions.
However, by \eqref{sloca} we can choose an integer $s\geqslant 0$ such that all products 
$(z_1 -z_2)^s\ts U_{i_1\ts j_1} \ldots U_{i_k\ts j_k}  $ modulo \eqref{modulo97}
are well-defined 
and such that
we have
$$
(z_1 -z_2)^s\ts T_{i_1\ts j_1} \ldots T_{i_k\ts j_k}  
=
(z_1 -z_2)^s\ts U_{i_1\ts j_1} \ldots U_{i_k\ts j_k} \mod \text{\eqref{modulo97}} .
$$

Let us turn to the second summand, $Y(z_2)(1\ot Y(z_1))$ in  \eqref{locality}. 
By \eqref{ymap}, its action on the  expression $x_{[m]}^{+23}(v)x_{[n]}^{+14}(u)$ produces
\beq\label{sloc3}
x_{[m]}^{23}(z_2 +v)\ts x_{[n]}^{13}(z_1 +u)
=x_{[m]}^{23}(z_2 +v_1,\ldots ,z_2 +v_m)\ts x_{[n]}^{13}(z_1 +u_1,\ldots , z_1+u_n).
\eeq
We now multiply \eqref{sloc3} by a sufficiently large power $(z_1-z_2)^r$, where $r\geqslant s$, so that in the 
  given expression, when regarded  modulo \eqref{modulo97}, we can employ  commutation relation \eqref{com12} to move each term $x(z_2+v_q)$  to the right of all $x(z_1 +u_p)$. Thus we get 
\beq\label{sloc4}
\sum_{k=0}^{\min\left\{m,n\right\}}
\sum_{ \substack{i=\left\{(i_1,j_1),\ldots ,(i_k,j_k)\right\}\in I_k^{n,m}\\l_1<\ldots< l_{n+m-2k} \in i' }}
(z_1-z_2)^r\ts
U_{i_1\ts j_1} \ldots  U_{i_k\ts j_k} \ts x_{[n-k]}^{13}\ts x_{[m-k]}^{23},
\eeq
where   $x_{[n-k]}^{13}$ and $ x_{[m-k]}^{23}$, as given by \eqref{sloc5} and \eqref{sloc6}, are applied on the tensor factors \eqref{sloc7} of \eqref{asso5}, respectively, and  the whole expression   is   regarded    modulo  \eqref{modulo97}.
Finally, let $y\in \hp$ be arbitrary. Apply \eqref{sloc2} and \eqref{sloc4} on $y$ and, furthermore, multiply the former expression by $(z_1-z_2)^r$. The preceding discussion shows that 
these two expressions,   both being regarded modulo \eqref{modulo97},  coincide,  
so that the $\Sc$-locality follows.
\end{prf}

\subsection{Quantum vertex algebra \texorpdfstring{$\Vc_{\hh}(c)$}{Vh(c)}}\label{subsec21}

Consider the topologically free subalgebra of $\hp$ generated by all  $x_{ii}^{(-r)}$ for $i=1,\ldots ,N$ and $r=1,2,\ldots $ and $\vac$. By using \eqref{normaln} and \eqref{es} one easily checks that this algebra is closed under the actions of the vertex operator map \eqref{ymap} and the braiding map \eqref{es}. Therefore,  Theorem \ref{mainthm} implies that this is a quantum vertex subalgebra of $\hp$. We denote this quantum vertex algebra by $\Vc_{\hh}(c)$. Clearly, as a $\CC[[h]]$-module $\Vc_{\hh}(c)$ coincides with the quotient of the
$h$-adically complete   algebra of polynomials
$$\CC[x_{ii}^{(-r)}\,:\, i =1,\ldots ,N,\, r\geqslant 1][[h]]$$
over its $h$-adically complete ideal generated by the elements in \eqref{34els}.

We now follow the exposition in \cite[Chap. 6]{LL} to recall the construction of the Heisenberg vertex algebra. Let $\h$ be an abelian Lie algebra with generators $a_1,\ldots ,a_N$ and the defining relations
$
a_1+\ldots +a_N =0$.
It is equipped with the  nondegenerate invariant symmetric
bilinear form $\left<\cdot,\cdot\right>$ given by
\beq\label{theform}
\left<a_i,a_j\right>=\delta_{ij}-\textstyle\frac{1}{N}\quad\text{for all }i,j=1,\ldots ,N. 
\eeq
The corresponding affine Lie algebra $\hhat$ is defined on the complex vector space
\eqref{introabafa}
via Lie brackets give by
\eqref{introbrckts}.  We organize all $a(r)=a\ot t^r$ in the  power series
$$
a(z)=\sum_{r\in\ZZ} a(r)\ts z^{-r-1}\quad\text{for all }a\in\h.
$$
Also, we write
$$
a^+(z)=\sum_{r\leqslant -1} a(r)\ts z^{-r-1}.
$$

Consider the subalgebras
$\hhat_{(\pm)} =\h\ot t^{\mp 1}\CC[t^{\mp 1}]\subset \hhat$.
For any $c\in\CC$ let
$$
V_{\hhat }(c)=U(\hhat) \ot_{U(\hhat_{(-)}\oplus \h \oplus\CC C)} \CC_c,
$$
be the induced $\hhat$-module, 
where $\hhat_{(-)}$ and $\h$ act trivially on $\CC_c =\CC$ and $C$ acts as the scalar $c$. 
By the Poincar\'{e}--Birkhoff--Witt theorem $V_{\hhat }(c)$ is isomorphic, as a vector space, to
the underlying space of the universal enveloping algebra 
$U(\hhat_{(+)})$. We can regard $\h$ as a subspace of $V_{\hhat }(c)$ via the map
\beq\label{sbspc4}
\h \ni a\mapsto a(-1)\vac\in V_{\hhat }(c),
\eeq
where $\vac=1\in\CC\subset V_{\hhat }(c)$. 
Using defining relations \eqref{introbrckts}  one can prove by induction the following identity   describing the action of  
$a(z)\in \ndo V_{\hhat }(c)[[z,z^{-1}]]$ with $a\in\h$
on $V_{\hhat }(c)$:
\begin{align}
&a(z)\ts b_1^+(z_1)\ldots b_n^+(z_n)
= 
a^+(z)\ts b^+_1(z_1)\ldots b^+_n(z_n)\non\\
&\qquad+\sum_{k=1}^n \frac{c \left<a,b_k\right> }{(z-z_k)^2}\ts
b^+_1(z_1)\ldots b^+_{k-1}(z_{k-1})\ts b^+_{k+1}(z_{k+1})\ldots b^+_n(z_n)) \label{clll4}
\end{align}
for any $b_1,\ldots ,b_n\in\h$.
The space $V_{\hhat }(c)$ can be equipped with the vertex algebra structure; cf. \cite{FZ, Lian}: 
\begin{thm}\label{llthm}
For any $c\in\CC$ there exists a unique vertex algebra structure on 
$V_{\hhat }(c)$ such that
$$
Y(a,z)=a(z)\in\ndo V_{\hhat }(c)[[z,z^{-1}]]\quad \text{for all}\quad a\in\h.
$$
\end{thm}

We now discuss the classical limit   of the quantum vertex algebra  $\Vc_{\hh}(c)$.
First, note that both $\Vc_{\hh}(c)$ and $V_{\hhat }(c)\equiv U(\hhat_{(+)}) $ can be naturally regarded as commutative associative algebras.
 Denote by
$\bar{x}_{ii}^{(-r)}$   the image of $x_{ii}^{(-r)}$ in the quotient $\Vc_{\hh}(c)_0\coloneqq \Vc_{\hh}(c) / h \Vc_{\hh}(c) $.
The assignments
$$
\bar{x}_{ii}^{(-r)}\mapsto a_i (-r) 
$$
with $r=1,2,\ldots$ and $i=1,\ldots ,N$  define an isomorphism
\beq\label{asgn2}
\Vc_{\hh}(c)_0 \to V_{\hhat }(c) 
\eeq
of complex commutative associative algebras. In particular,  for all $i=1,\ldots ,N$ map   \eqref{asgn2} identifies  $\bar{x}_{ii}^{(-1)}$    with $a_i=a_i (-i)\vac \in \h \subset V_{\hhat }(c)$; recall \eqref{sbspc4}. As the classical limit of  quantum vertex algebra is  vertex algebra (see \cite{EK}), let $\bar{Y}=\bar{Y}(z)$ be the classical limit of \eqref{ymap}, i.e. the vertex operator map of the vertex algebra $\Vc_{\hh}(c)_0$.
By Theorem \ref{llthm} the vertex operator map of $V_{\hhat }(c)$ is uniquely determined by its action on the elements of $ \h$. Therefore, in order to show that   \eqref{asgn2} is the vertex algebra isomorphism, it is sufficient to check that it maps  
\beq\label{sti7}
\bar{Y}(\bar{x}_{ii}^{(-1)},z)\bar{x}_{j_1 j_1}^{(-r_1)}\ldots \bar{x}_{j_n j_n}^{(-r_n)}\quad\text{to}\quad
Y(a_i,z)a_{j_1}(-r_1)\ldots a_{j_n}(-r_n)
\eeq
for all $n\geqslant 0$, $i,j_1,\ldots, j_n=1,\ldots ,N$ and $r_1,\ldots ,r_n=1,2,\ldots .$

Note that by \eqref{esform} the classical limit of the matrix entry $e_{ii}\ot e_{jj}$ of $S(z)$ equals 
$$
\textstyle\frac{c}{z^2}\left(\frac{1}{N}-\delta_{ij}\right)
=-\frac{c}{z^2}\left<a_i,a_j\right>.
$$
Hence, by \eqref{ymap},   the classical limit of the matrix entries
$e_{ii}\ot e_{j_1 j_1} \ot\ldots \ot e_{j_n j_n} $ 
in
\eqref{xm}
is
equal to
\begin{align}
&\bar{Y}(\bar{x}_{ii}^{(-1)},z)\ts \bar{x}^+_{j_1 j_1}(z_1)\ldots \bar{x}^+_{j_n j_n}(z_n) = \bar{x}_{ii}^+ (z)\ts \bar{x}^+_{j_1 j_1}(z_1)\ldots \bar{x}^+_{j_n j_n}(z_n)\non\\
 &\qquad+\sum_{k=1}^n \frac{c\left<a_i,a_{j_k}\right> }{(z-z_k)^2}\ts
 \bar{x}^+_{j_1 j_1}(z_1)\ldots \bar{x}^+_{j_{k-1} j_{k-1}}(z_{k-1})\ts
\bar{x}^+_{j_{k+1} j_{k+1}}(z_{k+1})\ldots \bar{x}^+_{j_n j_n}(z_n).\label{clll5}
\end{align}
Finally, by comparing the coefficients in  \eqref{clll4} and \eqref{clll5} we find that the second term in \eqref{sti7} coincides with the image of the first term under the  map \eqref{asgn2}.

\begin{pro}
For any $c\in\CC$ the classical limit $\Vc_{\hh}(c)_0$ of the quantum vertex algebra $\Vc_{\hh}(c)$ is the Heisenberg vertex algebra $V_{\hhat }(c) $.
\end{pro}

\section{Quantum Heisenberg algebra}\label{sec3}
\numberwithin{equation}{section}

In Subsection \ref{subsec31}, we introduce certain algebras $\hh(C) $ and $\hh(C)_* $. Moreover, we construct examples of their modules, which we use in Subsection \ref{subsec32}  to prove the Poincar\'{e}--Birkhoff--Witt theorem for these algebras. Finally, in Subsection \ref{subsec33}, we 
establish an equivalence between $\Vc_{\hh}(c)$-modules and certain class of $\hh(C) $-modules.

\subsection{Algebras \texorpdfstring{$\hh(C) $}{H(C)} and \texorpdfstring{$\hh(C)_* $}{H(C)*} }\label{subsec31}

The series  $S(u,C)$, as defined by \eqref{sofu},   can be written in the form
$$
S(u,C)=\sum_{i,j,k,l=1}^N e_{ij}\ot e_{kl}\ts  s_{ijkl}(u,C)\quad\text{for some}\quad s_{ijkl}(u,C)\in \CC[C,u^{-1}][[h]].
$$
To simplify the notation we write $s_{ij}(u,C)\coloneqq s_{iijj}(u,C)$.  Note that $s_{ij}(u,C) =s_{ji}(u,C)$. 

 Define  the algebra $\hh(C)$    as the $h$-adically complete associative algebra over the commutative ring $\CC[[h]]$ generated by the central element $C$ and the elements $y_i^{(r)}$, where $i=1,\ldots ,N$ and $r\in\ZZ,$ subject to defining relations written in terms  of the generator series
\beq\label{yi}
y^i (u) =\sum_{r\in\ZZ} y_i^{(r)} \ts u^{-r-1}\quad \text{for}\quad  i=1,\ldots ,N.
\eeq
The relations are given by $Cy=yC$ for all $y\in \hh(C)$ along with
\begin{gather}
y^i (u)\ts y^j(v)+s_{ij}(u-v,C)= y^j(v)\ts y^i (u)+s_{ij}(v-u,C),\label{relh1}\\
y^1(u)+\ldots +y^N(u)=0\label{relh2},
\end{gather}
where $i,j=1,\ldots ,N$.
It shall be useful to arrange the  generators into the diagonal matrix
\beq\label{genmatrix}
y(u)=\sum_{i=1}^N e_{ii}\ot y^i (u) \in\ndo\CC^N\ot \hh(C)[[u^{\pm 1}]].
\eeq
Also, we shall use the notation 
$$
\diag \left(\sum_{i_1,\ldots,i_n=1}^N \sum_{j_1,\ldots ,j_n=1}^N
e_{i_1 j_1}\ot\ldots\ot e_{i_n j_n}\right)
=
\sum_{i_1,\ldots,i_n=1}^N  
e_{i_1 i_1}\ot\ldots\ot e_{i_n i_n}
$$
to extract the diagonal matrix entries. Write $\Sd(u,C)=\diag S(u,C) $. The defining relations  \eqref{relh1} and \eqref{relh2} can be equivalently written in the matrix form as
\begin{gather*}
y_1(u) \ts y_2 (v) +\Sd (u-v,C)=y_2 (v)\ts y_1(u) +\Sd (v-u,C)\Fand 
\tr\ts  y(u) =0.
\end{gather*}

Let $W$ be an $\hh(C)$-module. Then, $W$ is said to be of {\em level $c$} if the action of $C$ on $W$ is scalar multiplication by $c\in \CC$. Moreover, $W$ is said to be {\em restricted} if it is a topologically free $\CC[[h]]$-module and the action of \eqref{yi} on $W$ satisfies
$$
y^i (z) \in\om(W,W((z))_h)\quad\text{for all}\quad i=1,\ldots ,N.
$$
If $W$ is   restricted module of level $c$,  the matrix of generators $y(u)$, as given by \eqref{genmatrix}, can be regarded as element of 
$
\ndo\CC^N \ot \om(W,W((u))_h).
$
Now denote by $u$ the family of variables $(u_1,\ldots ,u_n)$.
Following \eqref{normaln} one can introduce the elements
\beq\label{mod02}
y_{[n]}(u)=y_{[n]}(u_1,\ldots ,u_n)\in  (\ndo\CC^N)^{\ot n} \ot \om(W,W((u_1,\ldots , u_n))_h)
\eeq
for any $n=1,2,\ldots  $
by
$$
y_{[n]}(u) = 
\sum_{k=0}^{\floor{n/2}}
\sum_{ \substack{i=\left\{(i_1,j_1),\ldots ,(i_k,j_k)\right\}\in I^n_k\\l_1<\ldots< l_{n-2k} \in i' }}
\Sd_{i_1\ts j_1} \ldots  \Sd_{i_k\ts j_k} \ts y_{l_1}(u_{l_1})\ldots y_{l_{n-2k}}(u_{l_{n-2k}}),
$$
 where $\Sd_{ij} =\Sd_{ij} (u_i -u_j,c)$. As with \eqref{recall-}, we introduce the   elements  
\beq\label{citiraj}
y_{[n]}(z+u)\in(\ndo\CC^N)^{\ot n}\ot\om(\hp,\hp((z))_h [[u_1,\ldots ,u_n]])
\eeq
by
$$
y_{[n]}(z+u) = y_{[n]}(z+u_1,\ldots , z+u_n)\coloneqq
y_{[n]}(z_1,\ldots , z_n)\big|_{z_1=z+u_1,\ldots ,z_n=z+u_n}.\big.
$$

Note that due to \eqref{com12} and \eqref{trcom12} the assignments
\beq\label{asgn5}
C\mapsto c\Fand 
y_{i}^{(-r)}\mapsto x_{ii}^{(-r)}, \quad\text{where}\quad i=1,\ldots , N\fand r\in\ZZ,
\eeq
 define a structure of   
$\hh(C) $-module of  level $c$ on $\Vc_{\hh}(c)$. By \eqref{esform}  the series   $s_{ij}(u,c)$ belong  to $ \CC[u^{-1}][[h]]$ for all indices $i$ and $j$, so \eqref{xm} implies  that $\Vc_{\hh}(c)$ is   restricted.
This construction    can be generalized as follows. First, observe that by \eqref{esform} we have $\rez_u S(u) =0$. Hence the coefficient of $u_0^{-1}$ in \eqref{xm} is zero, which means that the action of $y_i^{(0)}$  on   $\Vc_{\hh}(c)$ is trivial. Furthermore, defining relation \eqref{relh1} implies 
$$ y_i^{(0)} y_j^{( r)}=y_j^{( r)}y_i^{(0)}\quad\text{for all}\quad  i,j=1,\ldots , N,\, r\in\ZZ.$$
Finally, note that the form \eqref{theform} can be uniquely extended to a $\CC[[h]]$-bilinear form on $\h[[h]]$, which we shall again denote by $\left<\cdot,\cdot\right>$. By the previous discussion and Lemma \ref{lemma21}, also recall  \eqref{sglop}, we have

\begin{pro}\label{lemma22}
For any $c\in \CC$ and $\alpha\in \h[[h]]$ there exists a unique structure of restricted $\hh(C)$-module of level $c$ on  $\Vc_{\hh}(c)$ such that for all $n\geqslant 0$ we have
$$
y_0(u_0) \diag \big(x_{[n]}^+ (u)\big)
=\diag \Big( \left(x _0(u_0)  +u_0^{-1} a_0 (\alpha) \right)     x_{[n]}^+ (u)\Big),
$$
where $u=(u_1,\ldots ,u_n)$ is a family of variables and $a(\alpha)=\sum_{i=1}^N \left<a_i,\alpha\right> \ts e_{ii}  $ the diagonal matrix
applied on the first tensor factor of $\ndo\CC^N \ot(\ndo\CC^N)^{\ot n}\ot \Vc_{\hh}(c)$.
\end{pro}

We denote the $\hh(C)$-module constructed in Proposition \ref{lemma22} by  $\Vc_{\hh}(c,\alpha)$. Clearly,    $\Vc_{\hh}(c,0)$ coincides with the $\hh(C)$-module  structure on  $\Vc_{\hh}(c )$
established by  
\eqref{asgn5}.

Let us show that the modules $\Vc_{\hh}(c,\alpha)$ present an $h$-adic generalization of the well-known canonical realization of the affine Lie algebra $\hhat$; see, e.g., \cite[Prop. 6.3.4]{LL}.
Suppose that the elements  $\wa_1,\ldots ,\wa_{N-1}$ form an orthonormal basis of the complex space $\h$, with respect to the bilinear form \eqref{theform}, such that
for all $j=1,\ldots ,N-1$ we have
$$
\wa_j=\sum_{i=1}^j \nu_{ij} \ts a_i \quad\text{for some}\quad\nu_{ij}\in\CC.
$$
Hence  for all $k=1,\ldots ,N-1$ we have
$$ 
\textstyle\xpan_\CC\left\{\wa_1,\ldots \wa_{k}\right\}=\xpan_\CC\left\{a_1,\ldots a_{k}\right\} .
$$ 
We employ   scalars $\nu_{ij}$ to introduce the elements  
\beq\label{wy1}
\wy_j^{(r)}=\sum_{i=1}^j \nu_{ij}\ts  y_i^{(r)}\in \hh(C), \quad\text{where}\quad j=1, \ldots ,N-1,\, r\in\ZZ.
\eeq
 Clearly, the elements \eqref{wy1} along with $C$ generate the algebra $\hh(C)$. Denote this particular set of generators by $\mathcal{G}_{\hh(C)}$.
Next, we define the corresponding family of   generators for the commutative algebra $\Vc_{\hh}(c)$ by
$$
\wx_{jj}^{(-r)}=\sum_{i=1}^j \nu_{ij}\ts  x_{ii}^{(-r)}\quad\text{with}\quad j=1, \ldots ,N-1,\, r=1,2,\ldots .
$$

\begin{kor}\label{cor22}
Let $c\in \CC$ and $\alpha\in \h[[h]]$. The action of the generators \eqref{wy1} on 
$$\Vc_{\hh}(c,\alpha)\equiv\CC [\ts\wx_{jj}^{(-r)}\ts :\ts j=1, \ldots ,N-1,\, r=1,2,\ldots\ts ][[h]]$$ 
is of the form
\begin{align*}
\wy_j^{(-r)} =\wx_{jj}^{(-r)}\mod h,\qquad
\wy_j^{(r)} =rc\frac{\partial}{\partial \wx_{jj}^{(-r)}}\mod h,\qquad
\wy_j^{(0)} = \left<\wa_j,\alpha\right>
\end{align*}
for all $j=1, \ldots ,N-1$ and $r=1,2,\ldots.$
\end{kor}

The algebra $\hh(C)$ may be regarded as an $h$-adic deformation of the universal enveloping algebra $U(\hat{\mathfrak{h}})$ of the commutative affine Lie algebra $\hat{\mathfrak{h}}$; see, in particular, the Poincar\'{e}--Birkhoff--Witt theorem for $\hh(C)$ below. Motivated by the classical theory, one can introduce the corresponding deformed Heisenberg algebra $\hh(C)_*$ in parallel with the definition of $\hh(C)$. More specifically,  $\hh(C)_*$ is the $h$-adically complete associative algebra over the commutative ring $\CC[[h]]$ generated by the central element $C$ and the elements $y_i^{(r)}$, where $i=1,\ldots ,N$ and $r\in\ZZ\setminus\left\{0\right\},$ subject to defining relations written in terms  of the generator series
$$
y^i (u) =\sum_{r\in\ZZ\setminus\left\{0\right\}} y_i^{(r)} \ts u^{-r-1}\quad \text{for} \quad i=1,\ldots ,N.
$$
Relations are given by $Cy=yC$ for all $y\in \hh(C)_*$ along with
\begin{gather*}
y^i (u)\ts y^j(v)+s_{ij}(u-v,C)= y^j(v)\ts y^i (u)+s_{ij}(v-u,C), \\
y^1(u)+\ldots +y^N(u)=0 ,
\end{gather*}
where $i,j=1,\ldots ,N$. 
The notion of restricted module of level $c$ for  $\hh(C)_*$ can be introduced in parallel with the corresponding definition for $\hh(C)$.
Furthermore, the suitable modifications of Proposition \ref{lemma22} and Corollary \ref{cor22} hold for the Heisenberg algebra as well.

\subsection{Poincar\'{e}--Birkhoff--Witt theorem}\label{subsec32}

For $c\in\CC$ denote by $\hh(c)$ the quotient of the algebra $\hh(C)$ over its $h$-adically closed ideal generated by $C-c$. It is clear from the defining relations \eqref{relh1} and \eqref{relh2} that the algebra $\hh(0)$ is isomorphic to the $h$-adically completed polynomial algebra in variables $\wy_j^{(r)}$, where $j=1, \ldots ,N-1$ and $r\in\ZZ$. Therefore, with  $\hh(0)$ being a quotient of $\hh(C)$, the ordered monomials in elements of $\mathcal{G}_{\hh(C)}$, with respect to any linear ordering, form a linearly independent subset in  $\hh(C)$.

Let us introduce a  linear ordering  on the set of generators $\mathcal{G}_{\hh(C)}$. Set 
\beq\label{order1}
C\prec \wy_j^{(r)}\qquad \text{for all}\qquad j=1, \ldots ,N-1\fand r\in\ZZ.
\eeq
Next, define 
\beq\label{order2}
\wy_i^{(s)}\prec \wy_j^{(r)}\qquad\text{if}\qquad s<r\quad\text{or}\quad s=r\text{ and }i<j.
\eeq
Let $\mathcal{M}_{\hh(C)}$ be the family of all increasing  monomials in generators, 
i.e. the monomials  of the form
$m_n\cdots m_1$, where   $n\geqslant 0$ and $m_1,\ldots, m_n\in \mathcal{G}_{\hh(C)}$, such that $m_n\prec\ldots\prec m_1$. 
We extend the ordering from $\mathcal{G}_{\hh(C)}$ to $\mathcal{M}_{\hh(C)}$ as follows.
For any two distinct monomials $m=m_n\cdots m_1$ and $m'=m'_{n'}\cdots m'_{1}$ in  
$\mathcal{M}_{\hh(C)}$
we write $m\prec m'$ if there exists an index $0\leqslant u\leqslant n+1,n'+1$ such that
$
m_1 =m'_1,\,\ldots,\, m_{u-1} =m'_{u-1}
$
and such that one of the following conditions holds:
$$n<u\leqslant n'\qquad\text{or}\qquad
 u \leqslant n,n'\fand m_u\prec m'_{u}.
$$
One easily checks that this defines a linear ordering on $\mathcal{M}_{\hh(C)}$.

The following theorem provides a topological basis of $\hh(C)$ in the $h$-adic sense.

\begin{thm}
The set $\mathcal{M}_{\hh(C)}$ forms a topological basis of $\hh(C)$.
\end{thm}

\begin{prf}
Note that defining relations \eqref{relh1} for the algebra $\hh(C)$, along with \eqref{wy1}, imply
\beq\label{pbw1}
\wy_i^{(r)}\wy_j^{(r)}=\wy_j^{(r)}\wy_i^{(r)}
\quad
\text{for all}\quad i,j=1,\ldots ,N-1,\, r\in\ZZ.
\eeq
Consider the family $F$ of all increasing monomials in  $C$ and
$y_j^{(r)}$, $j=1,\ldots ,N-1$, $r\in\ZZ$, where the ordering   is defined analogously to \eqref{order1} and \eqref{order2}, i.e. we have
$$C\prec y_j^{(r)}\qquad\Fand\qquad
y_i^{(s)}\prec y_j^{(r)}\qquad\text{if}\qquad  s<r \quad\text{or}\quad s=r\text{ and }i<j.$$
Here we assume that $F$ also contains the empty monomial, that is the unit $1$.
By using defining relations \eqref{relh1} and \eqref{relh2} one easily checks  that the $\CC[[h]]$-span of $F$ forms an $h$-adically dense $\CC[[h]]$-submodule of $\hh(C)$. However, by employing \eqref{wy1} and \eqref{pbw1}, we can express every monomial in $F$ as a $\CC$-linear combination of some elements of $\mathcal{M}_{\hh(C)}$. Thus, we conclude that the $\CC[[h]]$-span of $\mathcal{M}_{\hh(C)}$ is $h$-adically dense $\CC[[h]]$-submodule of $\hh(C)$ as well. 
Hence it remains to prove that $\mathcal{M}_{\hh(C)}$ is linearly independent over $\CC[[h]]$. 

Assume that in $\hh(C)$ we have a relation of the form 
\beq \label{linnez1}
\sum_{i \in I} p_i(C)m^i=0,
\eeq
where $I$ is a finite nonempty set, $p_i(x) \in \CC[[h]][x] $ are nonzero polynomials and $m^i \in  \mathcal{M}_{\hh(C)}$ are distinct monomials of the form $m^i=m_i^{+}m_i^{-}$, where $m_i^{-}$ (respectively $m_i^{+}$) are monomials in $\wy_j^{(r)}$ for $j=1,\ldots ,N-1$ and $r >0 $ (respectively $r <0$). 

\noindent  \hypertarget{ln-1}{(1)} Suppose that all $p_i(x)$ belong to  $\CC[x] $.
Let $i_0 \in I$ be such that for all $i \in I$ we have $m_{i_0}^{-} \prec m_{i}^{-}$ or $m_{i_0}^{-} = m_{i}^{-}$.  Denote by $J \subset I$   the set of all indices $j$ such that $m_{j}^{-}= m_{i_0}^{-}$.  Let $c \in \CC\hspace{-1pt}\hspace{-1pt} \setminus \hspace{-1pt}\hspace{-1pt}\{0\}$ be such that $p_j(c) \neq 0$ for all $j \in J$. Due to Corollary \ref{cor22}, we can choose an element $v \in \Vc_{\hh}(c,0)$ such that
\beq\label{linnez7}
m_{j}^{-} v \in\CC\hspace{-1pt}\hspace{-1pt}\setminus\hspace{-1pt}\hspace{-1pt}\left\{0\right\}+h \Vc_{\hh}(c,0)\text{ for all }j \in J
\fand  
m_{i}^{-} v \in h \Vc_{\hh}(c,0) \text{ for all } i \in I\hspace{-1pt}\setminus \hspace{-1pt} J.
\eeq
  If we act with 
	$$
	\sum_{i \in I} p_i(C)m^i=\sum_{j \in J} p_j(C)m_j^{+}m_j^{-} + \sum_{i \in I\setminus J} p_i(C)m_i^{+}m_i^{-}
	$$ on   $v$, from   \eqref{linnez7} follows the relation
\beq \label{linnez2}
\sum_{j \in J} p_j(c)(m_j^{-}v )m_j^{+} =0\mod h.
\eeq  
Note that all $p_j(c)(m_j^{-}v )$ belong to $\CC\hspace{-1pt}\hspace{-1pt}\setminus\hspace{-1pt}\hspace{-1pt}\left\{0\right\}$ modulo $h$.
 Since all    $m_j^{+}$ for $j \in J$  are mutually distinct ordered monomials   in generators of $\Vc_{\hh}(c,0)$, 
they are linearly independent, so that
 \eqref{linnez2} produces a contradiction.

\noindent \hypertarget{ln-2}{(2)} In   general case, i.e. for $p_i(x) \in \CC[[h]][x] $, choose a minimal integer $m\geqslant 0$ so that there exists a nonempty subset of indices $K\subset I$ such that
\begin{align*}
&p_k(x)\in h^m\CC[[h]][x]\hspace{-1pt}\hspace{-1pt}\setminus\hspace{-1pt}\hspace{-1pt}\left\{0\right\}+h^{m+1}\CC[[h]][x]\quad\text{for all}\quad k\in K,\\
&p_i(x)\in h^{m+1}\CC[[h]][x]\quad\text{for all}\quad i \in I\hspace{-1pt}\setminus\hspace{-1pt}  K.
\end{align*}
We   now suitably adapt the arguments from (\hyperlink{ln-1}{1}). Let $j_0 \in K$ be such that for all $k \in K$ we have $m_{j_0}^{-} \prec m_{k}^{-}$ or $m_{j_0}^{-} = m_{k}^{-}$.  Denote by $J \subset K$   the set of all indices $j$ such that $m_{j}^{-}= m_{j_0}^{-}$.
Choose
    $c \in \CC\hspace{-1pt}\hspace{-1pt} \setminus \hspace{-1pt}\hspace{-1pt}\{0\}$   such that $h^{-m}p_j(c) \in \CC \hspace{-1pt} \setminus \hspace{-1pt}\left\{0\right\}+h\CC[[h]]$ for all $j \in J$. As before, due to  Corollary \ref{cor22},  we can  choose  $v \in \Vc_{\hh}(c,0)$ such that
 \eqref{linnez7} holds. By applying \eqref{linnez1} on $v$ and then multiplying the equality  by $h^{-m}$ we get
\beq \label{linnez8}
\sum_{j \in J} (h^{-m}p_j(c))(m_j^{-}v )m_j^{+} =0\mod h.
\eeq  
Since   $(h^{-m}p_j(c))(m_j^{-}v )$ belongs to $\CC\hspace{-1pt}\hspace{-1pt}\setminus\hspace{-1pt}\hspace{-1pt}\left\{0\right\}$ modulo $h$ for all $j\in J$, \eqref{linnez8} leads to  contradiction as in (\hyperlink{ln-1}{1}).

So far we have proved that all  monomials in $\mathcal{M}_{\hh(C)}$ which do not contain elements $\wy_j^{(0)}$ are linearly independent over $\CC[[h]]$.
Assume that in \eqref{linnez1} all $m^i$ are mutually distinct   elements of $\mathcal{M}_{\hh(C)}$ of the form $m^i=m_i^{+}m_i^{0}m_i^{-}$, where $m_i^{\pm}$ are as before and $m_i^{0}$ are ordered monomials in $\wy_j^{(0)}$ for $j=1,\ldots ,N-1$. 
Such equality can be then written as
$$
\sum_{j \in J} q_j(C,\wy_1^{(0)},\ldots,\wy_{N-1}^{(0)} ) m_j^{+}m_j^{-}=0
$$
for some $J\subseteq I$ and nonzero polynomials $q_j(x_1,\ldots,x_N)\in\CC[[h]][x_1,\ldots ,x_N]$, so that all monomials $m_j^{+}m_j^{-}$ with $j\in J$ are mutually distinct. Indeed,  polynomials $q_j$ are found by
$$
q_j(C,\wy_1^{(0)},\ldots,\wy_{N-1}^{(0)} )=\sum_{ i\in I,\,m_i^{\pm} =m_j^{\pm} } p_i(C) m_i^0.
$$
However, we can choose $c\in \CC\hspace{-1pt}  \setminus  \hspace{-1pt}\{0\}$ and $\alpha\in\h$ so that the action of $q_j(C,\wy_1^{(0)},\ldots,\wy_{N-1}^{(0)} )$ on $\Vc_{\hh}(c,\alpha)$ is nonzero for all $j\in J$; recall Corollary \ref{cor22}.  Thus, the linear independence can be again established by  arguing as in (\hyperlink{ln-1}{1}) and (\hyperlink{ln-2}{2}).
\end{prf}

Denote by $\mathcal{M}_{\hh(C)_*}$ the family of all increasing monomials, with respect to the linear order $\prec$ defined by \eqref{order1} and \eqref{order2}, in generators
$C$ and $\wy_j^{(r)}$, where $j=1,\ldots , N-1$ and $r\in\ZZ\hspace{-1pt}\setminus\hspace{-1pt}\left\{0\right\}$,  of $\hh(C)_*$. The following corollary is clear.

\begin{kor}
The set $\mathcal{M}_{\hh(C)_*}$ forms a topological basis of $\hh(C)_*$.
\end{kor}

\subsection{Equivalence of  \texorpdfstring{$\hh(C) $}{H(C)}-modules and  \texorpdfstring{$\Vc_{\hh}(c)$}{Vh(c)}-modules}\label{subsec33}

The notion of module for   quantum vertex algebra   was   introduced by Li in parallel with the notion of vertex algebra module; see \cite[Def. 2.23]{Li}. By Lemma \ref{usefullemma} below, it coincides with the notion of module given by the following definition.

\begin{defn}\label{qvoamodule}
 Let $(V,Y,\vac ,\Sc)$  be a quantum vertex algebra. A {\em   $V$-module} is a pair $(W,Y_W)$, where $W$ is a topologically free $\CC[[h]]$-module
and $Y_W( z)$ a $\CC[[h]]$-module map
\begin{align*}
Y_W(z)\colon V\ot W&\to W((z))_h\\
v\ot w&\mapsto Y_W(z)(v\ot w)=Y_W(v,z)w=\sum_{r\in\mathbb{Z}} v_r w\ts z^{-r-1}
\end{align*}
which satisfies $Y_W(\vac,z)w=w$ for all $ w\in W$ and the {\em $\Sc$-Jacobi identity}
\begin{align}
&z_0^{-1}\delta\left(\frac{z_1 -z_2}{z_0}\right) Y_W(z_1)(1\ot Y_W(z_2))(u\ot v\ot w)\non\\
&\qquad-z_0^{-1}\delta\left(\frac{z_2-z_1}{-z_0}\right) Y_W(z_2)(1\ot Y_W(z_1)) \left(\Sc(-z_0)(v\ot u)\ot w\right)\non\\
&\qquad\qquad=z_2^{-1}\delta\left(\frac{z_1 -z_0}{z_2}\right)Y_W(Y(u,z_0)v,z_2 )w\quad\text{for all}\quad u,v \in V\text{ and }w\in W.\label{wjacobi}
\end{align}
 Let  $W_1$ be a topologically free $\CC[[h]]$-submodule of $W$. A pair $(W_1,Y_{W_1})$ is said to be a {\em $V$-submodule} of $W$ if $v_r w\in W_1$ for all $v\in V, w\in W_1$ and $r\in\ZZ$, where by $Y_{W_1}$ we denote  the restriction and corestriction of $Y_W$,
$$Y_{W_1}(z)=Y_W (z)\big|_{V\ot W_1}^{W_1}   \big. \,\colon\, V\ot W_1\,\to\, W_1 ((z))[[h]].$$
\end{defn}

The next  lemma  is well-known and follows by an argument similar to \cite[Rem. 2.16]{Li}; cf. also \cite[Lemma 1.3]{c11}. It is a quantum vertex algebra  analogue of \cite[Thm. 4.4.5]{LL}.

\begin{lem}\label{usefullemma}
Let $(V,Y,\vac,\Sc)$ be a quantum vertex algebra. Suppose  $W$ is a topologically free $\CC[[h]]$-module such that there exists a $\CC[[h]]$-module map
\begin{align*}
Y_W(z)\colon V\ot W&\to W((z))_h\\
v\ot w&\mapsto Y_W(z)(v\ot w)=Y_W(v,z)w=\sum_{r\in\mathbb{Z}} v_r w z^{-r-1}
\end{align*}
which satisfies $Y_W(\vac,z)w=w$  for all $w\in W$ and the {\em weak associativity}:
for any $u,v\in V$, $w\in W$ and $n\in\mathbb{Z}_{\geqslant 0}$
there exists $s\in\mathbb{Z}_{\geqslant 0}$
such that
\begin{align}
&(z_0 +z_2)^s\ts Y_W (u,z_0 +z_2)Y_W(v,z_2)\ts w\non\\
&\qquad - (z_0 +z_2)^s\ts Y_W\big(Y(u,z_0)v,z_2\big)\ts w
\in h^n W[[z_0^{\pm 1},z_2^{\pm 1}]].\label{associativityw}
\end{align}
Then $(W,Y_W)$ is a $V$-module.
In particular, 
it possesses the $\mathcal{S}$-{\em locality} property:
for any $u,v\in V$ and $n\in\mathbb{Z}_{\geqslant 0}$ there exists
$s\in\mathbb{Z}_{\geqslant 0}$ such that
\begin{align}
&(z_1-z_2)^{s}\ts Y_W(z_1)\big(1\otimes Y_W(z_2)\big)\big(\mathcal{S}(z_1 -z_2)(u\otimes v)\otimes w\big)
\nonumber\\
&-(z_1-z_2)^{s}\ts Y_W(z_2)\big(1\otimes Y_W(z_1)\big)(v\otimes u\otimes w)
\in h^n W[[z_1^{\pm 1},z_2^{\pm 1}]]\quad\text{for all }w\in W.\label{localityw}
\end{align}
\end{lem}

The goal of this section is to establish an equivalence between $\hh(c) $-modules and  $\Vc_{\hh}(c)$-modules. The next theorem is our first result in this direction.

\begin{thm}\label{mainthm2}
Let $W$ be an (irreducible) restricted $\hh(C) $-module of level $c\in\CC$. There exists a unique structure of (irreducible) $\Vc_{\hh}(c)$-module on $W$ given by  
\beq\label{mod01}
Y_W(\diag x^+_{[n]}(u_1,\ldots ,u_n),z)
=y_{[n]}(z+u_1,\ldots ,z+u_n)\quad\text{with}\quad n\geqslant 0.
\eeq
\end{thm}

\begin{prf}
Let $W$ be a restricted $\hh(C) $-module of level $c$.
First of all, we have to check that \eqref{mod01}, along with $Y_W(\vac, z)=1_W$, defines a $\CC[[h]]$-module map on $\Vc_{\hh}(c)$. It is sufficient to show that the ideals of relations
$$
 \big[x_{ii}^{(-r)},x_{jj}^{(-s)}\big]=0\fand
x_{11}^{(-r)}+\ldots + x_{NN}^{(-r)}=0,\quad\text{where}\quad i,j=1,\ldots ,N,\,r,s\geqslant 1
$$
are mapped to itself; recall Subsection \ref{subsec21}. This can verified by using the defining relations \eqref{relh1} and \eqref{relh2} for $\hh(C) $ and arguing as in the corresponding part of the proof of Theorem \ref{mainthm}; see Subsection \ref{appendix}. Furthermore, the map $Y_W(\cdot ,z)$ is uniquely determined by \eqref{mod01} as all monomials $x_{i_1 i_1}^{(-r_1)}\cdots x_{i_m i_m}^{(-r_m)}$ along with $\vac$ span an $h$-adically dense $\CC[[h]]$-submodule of $\Vc_{\hh}(c)$. Since $W$ is restricted, the image of $\Vc_{\hh}(c)$ under the map 
$v\mapsto Y_W(v,z)$ belongs to $\om (W,W((z))_h)$ due to \eqref{citiraj}.
Hence, in order to establish a structure of $\Vc_{\hh}(c)$-module on $W$ via \eqref{mod01}, it is sufficient to verify the weak associativity \eqref{associativityw}; recall Lemma \ref{usefullemma}. However, given the fact that the expressions $x_{[n]}(u)$ and $y_{[n]}(u)$ are defined by analogous formulae, \eqref{normaln} and \eqref{mod02} respectively, the   weak associativity can be proved by suitably adapting the arguments from the proof of Lemma \ref{assocl}. Thus we conclude that \eqref{mod01} defines a structure of $\Vc_{\hh}(c)$-module on $W$.

Finally, suppose that $W$ is an irreducible $\hh(C) $-module of level $c$. Let $W_1\subseteq W$ be a $\Vc_{\hh}(c)$-submodule of $W$. By employing \eqref{mod01} we find 
$$
y^i (z )w=Y_W (  x_{ii}^{(-1)},z)w\in W_1[[z^{\pm 1}]]\quad\text{for all }w\in W_1\text{ and }i=1,\ldots ,N,
$$
which implies $\hh(C)W_1\subseteq W_1$.
Thus,   $W_1$ is equal to $0$ or $W$, so that  $W$ is   irreducible as a  $\Vc_{\hh}(c)$-module as well.
\end{prf}

Note that by Theorem \ref{mainthm2} the $\hh(C)$-modules $\Vc_{\hh}(c,\alpha)$,
constructed in
Proposition \ref{lemma22}, can be equipped by the structure of module for the quantum vertex algebra $\Vc_{\hh}(c)$. 
We now show that the converse of Theorem \ref{mainthm2} also holds. Our proof relies on the $\Sc$-Jacobi identity   for quantum vertex algebra modules.

\begin{thm}\label{mainthm3}
Let $W$ be an (irreducible)   $\Vc_{\hh}(c)$-module for some $c\in\CC$. There exists a unique structure of (irreducible) restricted $\hh(C) $-module of level $c$ on $W$ such that  
\beq\label{mod03}
y^i (z ) =Y_W (  x_{ii}^{(-1)},z) \quad\text{for all} \quad i=1,\ldots ,N .
\eeq
\end{thm}

\begin{prf}
Let $W$ be a $\Vc_{\hh}(c)$-module. The  map $Y_W(\cdot, z)$ satisfies the   $\Sc$-Jacobi identity \eqref{wjacobi},
\begin{align}
&z_0^{-1}\delta\left(\frac{z_1 -z_2}{z_0}\right) Y_W(z_1)(1\ot Y_W(z_2)) \non\\
&\qquad-z_0^{-1}\delta\left(\frac{z_2-z_1}{-z_0}\right) Y_W(z_2)(1\ot Y_W(z_1)) \left(\Sc(-z_0)P\ot 1 \right)\non\\
   =&\,z_2^{-1}\delta\left(\frac{z_1 -z_0}{z_2}\right)Y_W(z_2) (Y(z_0)\ot 1  )  \label{wjacobi2}
\end{align}
on $\Vc_{\hh}(c)\ot \Vc_{\hh}(c)\ot W$. First, by applying the identity on
$$
\diag\left(x_{13}^+(0)\ot x_{24}^+(0)\ot w\right)=\sum_{i,j=1}^N e_{ii}\ot e_{jj} \ot x_{ii}^{(-1)}\ot x_{jj}^{(-1)}\ot w,
$$
where $w\in W$ is arbitrary,
and using the definitions of the corresponding maps, recall \eqref{es}, \eqref{ymap} and \eqref{mod03},    we get
\begin{align}
&z_0^{-1}\delta\left(\frac{z_1 -z_2}{z_0}\right) Y_W(\diag x_{13}^+(0),z_1)  Y_W(\diag x_{23}^+(0),z_2)w \non\\
&\qquad-z_0^{-1}\delta\left(\frac{z_2-z_1}{-z_0}\right)  \left(Y_W(\diag x_{23}^+(0),z_2)  Y_W(\diag x_{13}^+(0),z_1)w+ \bar{T}_{12}(-z_0,c) \ot w\right)\non\\
= & \,z_2^{-1}\delta\left(\frac{z_1 -z_0}{z_2}\right)\left(Y_W( \diag(x_{13}^+(z_0)x_{23}^+(0)),z_2 )w -\bar{S}_{12}(z_0,c)\ot w\right) , \label{wjacobi3}
\end{align}
where $\bar{T}(z,c)=\bar{S}(z,c)-\bar{S}(-z,c)$; recall \eqref{te}.  
Next, by the property of the delta function we have
\begin{align}
 \left(z_0^{-1}\delta\left(\frac{z_1 -z_2}{z_0}\right)    -z_0^{-1}\delta\left(\frac{z_2-z_1}{-z_0}\right)\right)  \bar{S}_{12}(z_0,c) \ot w 
=  z_2^{-1}\delta\left(\frac{z_1 -z_0}{z_2}\right)\bar{S}_{12}(z_0,c) \ot w . \label{wjacobi4}
\end{align}
Finally, by adding equalities \eqref{wjacobi3} and \eqref{wjacobi4} and employing \eqref{te} we obtain
\begin{align}
&z_0^{-1}\delta\left(\frac{z_1 -z_2}{z_0}\right) \left(Y_W(\diag x_{13}^+(0),z_1)  Y_W(\diag x_{23}^+(0),z_2)+ \bar{S}_{12}(z_0,c)\right)w \non\\
&\qquad-z_0^{-1}\delta\left(\frac{z_2-z_1}{-z_0}\right)  \left(Y_W(\diag x_{23}^+(0),z_2)  Y_W(\diag x_{13}^+(0),z_1) + \bar{S}_{12}(-z_0,c)  \right)w\non\\
= & \,z_2^{-1}\delta\left(\frac{z_1 -z_0}{z_2}\right) Y_W( \diag(x_{13}^+(z_0)x_{23}^+(0)),z_2 )w   . \non
\end{align}
Taking the residue with respect to the variable $z_0$ produces the identity
\begin{align}
&   Y_W(\diag x_{13}^+(0),z_1)  Y_W(\diag x_{23}^+(0),z_2)+ \bar{S}_{12}(z_1-z_2,c) \non \\
&\qquad-    Y_W(\diag x_{23}^+(0),z_2)  Y_W(\diag x_{13}^+(0),z_1)w+ \bar{S}_{12}(z_2-z_1,c)   
= 0    \label{wjacobi5}
\end{align}
which holds when applied on $w$. With $w$ being an arbitrary element of $W$, we conclude by \eqref{wjacobi5}  that the map   \eqref{mod03} satisfies defining relation \eqref{relh1}. As for the other defining relation  \eqref{relh2},  it is also satisfied by  the map \eqref{mod03} due to \eqref{trcom12}. In addition,  Definition \ref{qvoamodule} implies that the image of   $Y_W(\cdot,z)$ belongs to $\om (W,W((z))_h)$. Therefore, formula \eqref{mod03} defines a structure of restricted $\hh(C) $-module of level $c$ on $W$, as required.

As for the second assertion of the theorem, suppose that $W$ is an irreducible $\Vc_{\hh}(c)$-module. Let $W_1\subseteq W$ be a  $\hh(C)$-submodule of $W$. 
By using \eqref{mod01} we get
\beq\label{nice2}
Y_W(\diag x^+_{[n]}(u ),z)w
=y_{[n]}(z+u )w
\in W_1[[z^{\pm 1},u_1,\ldots ,u_n]]
\eeq
 for all    $n\geqslant 1$ and  $w\in W_1$, where $u=(u_1,\ldots ,u_n)$; cf. Remark \ref{nice} below.
This implies  that $  W_1$ is a $\Vc_{\hh}(c)$-submodule, so we conclude that 
$W_1$ equals $0$ or $W_1$. Hence $W $ is an irreducible $\hh(C)$-module,  as required.
\end{prf}

\begin{rem}\label{nice}
We should say that in the end of the proof of Theorem \ref{mainthm3} we     use  the fact that the $\Vc_{\hh}(c)$-module map $Y_W$ coming from  \eqref{mod03} takes the form \eqref{nice2}. This   is easily verified by using the so-called iterate formula,
\begin{align}
Y_W(Y(u,z_0)v,z_2 )w=  
&\ts\rez_{z_1}\left(z_0^{-1}\delta\left(\frac{z_1 -z_2}{z_0}\right) Y_W(z_1)(1\ot Y_W(z_2))(u\ot v\ot w)\right.\non\\
&\left.  -z_0^{-1}\delta\left(\frac{z_2-z_1}{-z_0}\right) Y_W(z_2)(1\ot Y_W(z_1)) \left(\Sc(-z_0)(v\ot u)\ot w\right)\right),\non 
\end{align}
 which is obtained by taking the residue $\rez_{z_1}$ of the $\Sc$-Jacobi identity \eqref{wjacobi}.
More specifically, due to the iterate formula, it is sufficient to check that the $h$-adic completion of the $\CC[[h]]$-module $S\subset \Vc_{\hh}(c)$, which contains the elements $x_{11}^{(-1)},\ldots, x_{NN}^{(-1)}$ and $\vac$ and satisfies $u_r v\in S$ for all $u,v\in S$ and $r\in \ZZ$, coincides with $\Vc_{\hh}(c)$. However, one can prove by induction over $n$ that all monomials   $x_{i_1\ts i_1}^{(-r_1)}\cdots x_{i_n\ts i_n}^{(-r_n)}$ with $i_1,\ldots ,i_n=1,\ldots, N$ and $r\geqslant 1$ belong to $S$, which   implies such conclusion. Indeed, for $n=1$ this follows by extracting the coefficients of $Y(x_{ii}^{(-1)},z) \vac$. Suppose that the statement holds for all monomials of length less than or equal to $n$. By \eqref{normaln} the coefficients of matrix entries of
\begin{align*}
&\diag x_{[n+1]}^+(z,v_1,\ldots,v_n)
\Fand\\
&\diag
x^{13}(z)x_{[n]}^{+23}(v_1,\ldots ,v_n)
=\diag Y(\ts x^{+13}(0),z)x_{[n]}^{+23}(v_1,\ldots ,v_n)
\end{align*}
coincide modulo monomials of length less than or equal to $n$, so the induction assumption implies that the statement holds for the monomials of length $n+1$ as well.
\end{rem}

\section*{Acknowledgement}
This work has been supported in part by Croatian Science Foundation under the project UIP-2019-04-8488. 
 The first  author is partially supported by the QuantiXLie Centre of Excellence, a project cofinanced by the Croatian Government and European Union through the European Regional Development Fund - the Competitiveness and Cohesion Operational Programme (Grant KK.01.1.1.01.0004).


\end{document}